\documentclass[12pt,oneside,reqno]{amsart}

\usepackage{txfonts}
\usepackage{bbm}

\usepackage{amsmath}
\usepackage{graphicx}
\usepackage{mathrsfs}
\usepackage{stmaryrd}
\usepackage{amsfonts}
\usepackage{enumerate,amsmath,amssymb,amsthm,color}
\numberwithin{equation}{section}

\newcommand{\be}{\begin{eqnarray}}
\newcommand{\ee}{\end{eqnarray}}
\newcommand{\ce}{\begin{eqnarray*}}
\newcommand{\de}{\end{eqnarray*}}
\newtheorem{theorem}{Theorem}[section]
\newtheorem{lemma}[theorem]{Lemma}
\newtheorem{remark}[theorem]{Remark}
\newtheorem{definition}[theorem]{Definition}
\newtheorem{proposition}[theorem]{Proposition}
\newtheorem{Examples}[theorem]{Example}
\newtheorem{corollary}[theorem]{Corollary}

\def\e{\mathrm{e}}

\def\p{\partial}

\def\[{{\Big[}}
\def\]{{\Big]}}
\def\<{{\langle}}
\def\>{{\rangle}}
\def\({{\Big(}}
\def\){{\Big)}}

\def\bx{{\mathbf{x}}}

\def\dif{{\mathord{{\rm d}}}}

\def\no{\nonumber}
\def\={&\!\!=\!\!&}
\def\bt{\begin{theorem}}
\def\et{\end{theorem}}
\def\bl{\begin{lemma}}
\def\el{\end{lemma}}
\def\br{\begin{remark}}
\def\er{\end{remark}}
\def\bd{\begin{definition}}
\def\ed{\end{definition}}
\def\bp{\begin{proposition}}
\def\ep{\end{proposition}}
\def\bc{\begin{corollary}}
\def\ec{\end{corollary}}
\def\bx{\begin{Examples}}
\def\ex{\end{Examples}}

\def\mE{{\mathbb E}}

\def\mI{{\mathbb I}}

\def\mN{{\mathbb N}}

\def\mP{{\mathbb P}}

\def\mR{{\mathbb R}}

\def\1{{\mathbf{1}}}

\def\sL{{\mathscr L}}

\def\geq{\geqslant}
\def\leq{\leqslant}

\allowdisplaybreaks

\def\be{\begin{equation}}
\def\ee{\end{equation}}

\begin{document}

\title{Stochastic flows for L\'evy processes with H\"{o}lder drifts}

\date{}
\author{{Zhen-Qing Chen}, \   {Renming Song} \  and  \ {Xicheng  Zhang}}

\address{Zhen-Qing Chen:
Department of Mathematics, University of Washington, Seattle, WA 98195, USA   \newline
Email: zqchen@uw.edu
 }
\address{Renming Song: Department of Mathematics, University of Illinois, Urbana, IL 61801, USA \newline
Email: rsong@math.uiuc.edu
}
\address{Xicheng Zhang:
School of Mathematics and Statistics, Wuhan University,
Wuhan, Hubei 430072, P.R.China \newline
Email: XichengZhang@gmail.com
 }

\thanks{The research of ZC is partially supported
by NSF grant DMS-1206276.
The research of RS is partially supported by a grant from the Simons Foundation (208236).
 The research of XZ is partially supported by NNSFC grant of China (Nos. 11271294, 11325105).}

\begin{abstract}
In this paper we study the following stochastic differential equation (SDE) in $\mR^d$:
$$
\dif X_t= \dif Z_t + b(t, X_t)\dif t, \quad  X_0=x,
$$
where $Z$ is a L\'evy process.
We show that for a large class of L\'evy processes ${Z}$ and H\"older continuous
drift $b$, the SDE above has a unique strong solution for every starting point $x\in\mR^d$.
Moreover, these strong solutions form a $C^1$-stochastic flow.
As a consequence, we show that, when ${Z}$ is
an $\alpha$-stable-type L\'evy
process with $\alpha\in (0, 2)$ and $b$ is
bounded and $\beta$-H\"older continuous
with $\beta\in (1- {\alpha}/{2},1)$, the SDE above has a unique strong solution.
When $\alpha \in (0, 1)$, this in particular
solves an open problem from Priola \cite{Pr1}.
Moreover, we obtain a Bismut type derivative formula
for $\nabla \mE_x f(X_t)$ when
 ${Z}$ is a subordinate Brownian motion.
To study the SDE above, we first study the following
nonlocal parabolic  equation with H\"older continuous $b$ and $f$:
$$
\p_t u+\sL u+b\cdot \nabla u+f=0,\quad u(1, \cdot )=0,
$$
where $\sL$ is the generator of the L\'evy process ${Z}$.

\bigskip

\noindent {{\bf AMS 2010 Mathematics Subject Classification:} Primary 60H10,  35K05;   Secondary 60H30,   47G20}

\noindent{{\bf Keywords and Phrases:}
SDE, supercritical, subcritical, stable process, subordinate Brownian motion,
gradient estimate, strong existence, pathwise uniqueness, Bismut formula, stochastic flow,
$C^1$-diffeomorphism}

\end{abstract}

\maketitle
\rm

\section{Introduction}

Consider the following stochastic differential equation (SDE) in $\mR^d$:
\begin{align}
\dif X_t= \dif Z_t + b (t, X_t)\dif t, \quad  X_0=x,\label{SDE1}
\end{align}
where
$b(t, x): [0,1]\times\mR^d\to\mR^d$ is a
bounded Borel function
and ${Z}$ is a L\'evy process in $\mR^d$.
When $d=1$, ${Z}$ is a Brownian motion and $b$ is a bounded Borel function on $\mR$,
Zvonkin \cite{Zv} proved that the above SDE admits a unique strong solution for every starting point $x$.
Zvonkin's result was extended to the multi-dimensional case by Veretennikov \cite{Ve}.
Since then, many people have made contributions to this problem
(see \cite{Kr-Ro, Fl-Gu-Pr, Fe-Fl, Zh0} and references therein).
However, when ${Z}$ is a pure jump L\'evy process, strong existence and pathwise uniqueness
of SDE (\ref{SDE1}) become quite involved for drift $b$ which is not Lipschitz continuous.
When $d=1$, $b(t, x)=b(x)$ and
${Z}$ is a symmetric $\alpha$-stable process in $\mR$ with $\alpha\in(0,1)$,
Tanaka, Tsuchiya and Watanabe \cite{Ta-Ts-Wa}
proved that pathwise uniqueness fails for \eqref{SDE1} even for bounded
$b\in C^\beta_b(\mR)$.
On the other hand, when $d=1$ and ${Z}$ is a symmetric $\alpha$-stable process in
$\mR$ with $\alpha\in[1,2)$, it is shown in \cite{Ta-Ts-Wa}
that  pathwise uniqueness holds for
\eqref{SDE1} for any bounded continuous $b(t, x)=b(x)$.
For $d\geq 2$, using Zvonkin's transform, Priola \cite{Pr} obtained pathwise
uniqueness for SDE (\ref{SDE1}) when ${Z}$ is a
non-degenerate symmetric (but possibly non-isotropic)
$\alpha$-stable process in $\mR^d$ with $\alpha\in[1,2)$  and
time-independent $b(t, x)=b(x) \in C^\beta_b(\mR^d)$
with $\beta\in(1- {\alpha}/{2},1)$.
Note that in this case, the infinitesimal generator corresponding to the
solution $X$ of \eqref{SDE1} is
$\sL^{(\alpha)}+b\cdot \nabla$. Here $\sL^{(\alpha)}$ is the infinitesimal generator
of the L\'evy process ${Z}$, which is a nonlocal operator of order $\alpha$.
When $\alpha>1$, $\sL^{(\alpha)}$
is the dominant term, which is called the subcritical case.
When $\alpha \in (0, 1)$, the  gradient $\nabla$ is of higher order than
the nonlocal operator $\sL^{(\alpha)}$
so the corresponding SDE \eqref{SDE1} is called  supercritical. The critical case corresponds to $\alpha=1$.
Priola's result was extended to drifts $b$ in some
fractional Sobolev spaces in the subcritical case in \cite{Zh1}
and to more general L\'evy processes in the subcritical case in \cite{Pr1}.
However, when $d\geq 2$,
$\alpha\in(0,1)$ and ${Z}$ is a symmetric $\alpha$-stable process in $\mR^d$,
even for time-independent H\"older continuous drift $b$,
pathwise uniqueness for SDE (\ref{SDE1}) was an open question until now;
see \cite[Remark 5.5]{Pr1}.
When ${Z}$ is a rotationally symmetric $\alpha$-stable process, SDE (\ref{SDE1})
is connected with the following nonlocal PDE:
$$
\p_t u+\Delta^{\alpha/2}u+ b \cdot \nabla u +f =0,
$$
where $\Delta^{\alpha/2}:=
-(-\Delta)^{\alpha/2}$ is the usual fractional Laplacian.
In order to solve SDE (\ref{SDE1}) driven by a rotationally symmetric stable process ${Z}$,
one needs to understand the above PDE better.
In this direction, Silvestre \cite{Si2}
obtained  the following a priori interior estimate:
$$
\|u\|_{L^\infty([0,1];
C^{\alpha+\beta}(B_1))}
\leq C\Big(\|u\|_{L^\infty([0,2]\times B_2)}+\|f\|_{L^\infty([0,2]; C^{\beta}(B_2))}\Big),
$$
where, for any $r>0$, $B_r$ stands for the open ball of radius $r$ centered at the origin,
provided $b\in L^\infty([0,2]; C^{\beta}(B_2))$ and $\alpha+\beta>1$. Such an estimate
suggests that one could solve
the supercritical
SDE (\ref{SDE1}) uniquely
when ${Z}$ is a rotationally  symmetric $\alpha$-stable process with $\alpha\in(0,1)$
and  $b\in C^\beta_b(\mR^d)$ with $\beta\in(1-{\alpha}/{2},1)$.
However, this is not an easy task since one needs
additional asymptotic estimates
in the time variable.
Furthermore, the approach of \cite{Si2} strongly depends on
realizing the fractional Laplacian in $\mR^d$ as
the boundary trace of an elliptic operator in upper half space of $\mR^{d+1}$.
Extending this approach to other nonlocal operators,
such as $\alpha$-stable-type operators, would be very hard if not impossible.

The goal of this paper is to establish strong existence and pathwise uniqueness
for SDE (\ref{SDE1}) with, possibly time-dependent, H\"older continuous drift $b$
for a large class of L\'evy processes including stable-type L\'evy processes.
We not only extend the main result of \cite{Pr1} in the subcritical case to more general
L\'evy processes and time-dependent drifts
 but also establish strong existence and pathwise uniqueness result in
the supercritical case for a large class of L\'evy processes
where the drift $b$ can be time-dependent.
We emphasize that the L\'evy process ${Z}$ in this paper can be non-symmetric
and may also have drift.
One of the main results
of this paper in particular solves the open problem raised in
\cite[Remark 5.5]{Pr1} where ${Z}$ is a symmetric $\alpha$-stable process with
$\alpha \in (0, 1)$.
Our approach is mainly probabilistic.

In this paper, we use ``:=" as a way of definition. For $a, b\in \mR$,
$a\wedge b:=\min\{a, b\}$, $a\vee b:=\max\{a, b\}$, and $a^+:=a\vee 0$.
We now describe the setup and
the main results of this paper.
Let $\sL_{\nu, \eta}$
be the infinitesimal generator of the L\'evy process ${Z}$, that is,
$$
\sL_{\nu, \eta} f(x)= \int_{\mR^d}\left( f(x+z)-f(x)- \1_{\{|z|\leq 1\}}z\cdot\nabla f(x)\right)
\nu(\dif z) +\eta\cdot\nabla f(x),
$$
where $\nu$ is the L\'evy measure of ${Z}$ and $\eta$ is a vector in $\mR^d$.
For any $\eta\in \mR^d$ and any L\'evy measure $\nu$, i.e.,
a measure on $\mR^d\setminus\{0\}$ with $\int(1\wedge |z|^2)\nu(\dif z)<\infty$,
we will use
$\{{T}^{\nu, \eta}_t; t\geq 0\}$ to denote the transition
semigroup of the L\'evy process ${Z}$ with infinitesimal generator  $\sL_{\nu, \eta}$, i.e.,
$$
{T}^{\nu, \eta}_tf(x):=\mE f(x+Z_t).
$$
Suppose that $\nu$ can be decomposed as
\begin{align}
\nu=\nu_0+\nu_1+\nu_2,\label{Dec}
\end{align}
where $\nu_1, \nu_2$ are two L\'evy measures, and $\nu_0$ is a
{\it finite signed} measure so that
$$
\nu_0+\nu_1\mbox{ is still a  L\'evy measure}.
$$
We make the following assumption about ${T}^{\nu_1, 0}_t$.
There exist $\alpha\in (0, 2)$,  $\bar\alpha, \delta\in(0,1]$ and $K_0>0$
so that the following  gradient estimates for
the semigroup
$\{T^{\nu_1, 0}_t; t\geq 0\}$
hold.
\medskip

\begin{enumerate}
\item  [{\bf \big(H$^{\alpha,\bar\alpha,\delta}_{\nu_1,K_0}$\big)}]
If $\alpha\in(0,1]$, then for any $x\in\mR^d$, $\beta\in [0, \bar\alpha]$ and
bounded Borel function $f$ satisfying
$$
|f(x+y)-f(x)|\leq \Lambda|y|^\beta  \quad \hbox{for all } y\in\mR^d,
$$
with some $\Lambda>0$, it holds that
\begin{align}
|\nabla{T}^{\nu_1, 0}_t f(x)|\leq K_0\Lambda t^{(\delta\beta-1)/\alpha}
\quad  \hbox{for all } t\in(0,1).  \label{Es13}
\end{align}
\item  [{\bf \big(H$^{\alpha}_{\nu_1,K_0}$\big)}]
If $\alpha\in(1,2)$, then for any bounded Borel function $f$, it holds that
\begin{align}
\|\nabla{T}^{\nu_1,0}_t f\|_\infty\leq K_0\|f\|_\infty t^{-1/\alpha}
\quad  \hbox{for all } t\in(0,1). \label{Es133}
\end{align}
 \end{enumerate}

\br\rm
The pointwise estimate \eqref{Es13} allows us to borrow the H\"older regularity of
the drift to compensate the time singularity,
which is crucial for the well-posedness of SDEs with H\"older drifts in the supercritical case.
Condition (\ref{Es133}) in the subcritical case is the same as  Hypothesis 1 of Priola \cite{Pr1}.
Moreover, the parameters $\bar\alpha$ and $\delta$ are mainly designed for Example \ref{E:4.4} below, and
in Examples \ref{EEU} and \ref{EEU1} below, $\bar\alpha$ and $\delta$
can all be chosen to be $1$.
\er

The first main result of this paper is the following

\bt\label{Main1}
Suppose that either {\bf \big(H$^{\alpha,\bar\alpha,\delta}_{\nu_1,K_0}$\big)} holds for
some $\alpha\in(0,1]$,  $\bar\alpha, \delta\in (0, 1]$ and $K_0>0$, or
{\bf \big(H$^{\alpha}_{\nu_1,K_0}$\big)} holds for some $\alpha\in(1,2)$ and $K_0>0$.
Let $\gamma \in (0, 1)$ be such that
 \begin{align}
\int_{|z|\leq 1}|z|^{2\gamma}\nu(\dif z)<\infty.\label{EB1}
\end{align}
Assume further that $\gamma + (1-\alpha)/\delta<\bar \alpha$
in the case $\alpha\in (0, 1]$.
If for some
$\beta\in(\gamma + (1-\alpha)/\delta, 1]$
in the case $\alpha\in(0,1]$, and for some
$\beta\in((\gamma + 1-\alpha)^+, 1]$
in the case $\alpha\in(1,2)$,
it holds that
\begin{align}\label{Ew22}
\sup_{t\in [0, 1]} \|b(t, \cdot)\|_\infty+\sup_{t\in[0,1]}\sup_{x\not=y\in\mR^d}
\frac{|b (t, x)-b (t, y)|}{|x-y|^\beta}<\infty,
\end{align}
then for every $x\in \mR^d$, there is a unique strong solution
$\{X_t(x); t\in[0,1]\}$
to SDE (\ref{SDE1}). Moreover, $\{X_t(x), t\in[0,1], x\in\mR^d\}$ forms a $C^1$-stochastic
diffeomorphism flow, and
for each $x\in\mR^d$, $t\mapsto \nabla X_t(x)$ is continuous, and
\begin{align}
\sup_{x\in\mR^d}\mE\left[\sup_{t\in[0,1]}|\nabla X_t(x)|^p\right]
\leq C_p <\infty \quad  \hbox{for every }  p\geq 1,\label{EQ9}
\end{align}
where $C_p$ only depends on
$p,d,\alpha,\beta,\gamma,\nu, K_0, \bar \alpha, \delta$
and the H\"older norm of $b$.
\et
\br

\rm By a suitable localization argument (cf. \cite{Zh1}), for the local uniqueness of
SDE \eqref{SDE1}, the global condition \eqref{Ew22} can be replaced with a local condition.
Moreover, although $t\mapsto X_t(x)$ is not continuous, since we are considering an additive noise, the conclusion that $t\mapsto\nabla X_t(x)$ is continuous is not surprising.
\er

Various examples of L\'evy processes satisfying the conditions
{\bf \big(H$^{\alpha,\bar\alpha,\delta}_{\nu_1,K_0}$\big)} with $\alpha \in (0, 1]$,
{\bf \big(H$^{\alpha}_{\nu_1,K_0}$\big)} with $\alpha \in (1, 2)$,
and \eqref{EB1}
(and hence the conclusion of Theorem \ref{Main1} holds for these L\'evy processes)
are given in Section 4.
To illustrate Theorem \ref{Main1}, here we only give the following corollary,
which is a direct consequence of these examples.

\begin{corollary} \label{C:1.3}

\begin{enumerate}[{\rm (i)}]
\item {\bf (Stable-type L\'evy process)}
 Let ${Z}$ be a L\'evy process in $\mR^d$
whose L\'evy measure $\nu$ has a density $\kappa(z)$.
Assume that for some
$0<\alpha_1\leq\alpha_2<2$,
$$
c_1|z|^{-d-\alpha_1}\leq \kappa(z)\leq c_2|z|^{-d-\alpha_2}
\ \ \hbox{for }\
0<|z|\leq 1.
$$
Assume that $\alpha_2<2\alpha_1$, and $b(t, x)$ is bounded and  $\beta$-H\"older continuous
in $x$  uniformly in $t\in [0, 1]$,  for some
$\beta \in (1+\alpha_2/2 -\alpha_1 , 1]$.
Then SDE \eqref{SDE1} has a unique strong solution
for every $x\in \mR^d$ and \eqref{EQ9} holds.

\smallskip

\item {\bf (Subordinate Brownian motion)} Let ${Z}$ be a subordinate Brownian motion in
$\mR^d$ with characteristic function $\Phi (z)$. Suppose that
there are $0<\alpha_1\leq\alpha_2<2$ such that
$$
C_1 |z|^{ \alpha_1}\leq \Phi (z) \leq C_2 |z|^{\alpha_2}\quad \hbox{for } |z| \geq 1.
$$
Assume that $\alpha_2<2\alpha_1$, and $b(t, x)$ is bounded and  $\beta$-H\"older continuous
in $x$ uniformly in $t\in [0, 1]$,  for some
$\beta \in (1+ \alpha_2/2 -\alpha_1 , 1]$.
Then SDE \eqref{SDE1} has a unique strong solution for every $x\in \mR^d$ and \eqref{EQ9} holds.

\smallskip

\item { \bf (Cylindrical stable process)}
Let $Z =(Z^1, \cdots, Z^k)$,
where $Z^j$, $1\leq j\leq k$, are independent $d_j$-dimensional
rotationally symmetric $\alpha_j$-stable processes, respectively,
with $\alpha_j\in (0, 2)$ and $d_j\geq 1$.
Let $\alpha:=\min_{1\leq j\leq k}\alpha_j$ and $\alpha_{\mathrm{max}}:=\max_{1\leq j\leq k}\alpha_j$.
Suppose that
\begin{equation}\label{e:1.9}
\hbox{either} \quad \alpha > 1 \quad \hbox{or} \quad
\alpha \in (0, 1] \ \hbox{ and } \ \alpha_{\mathrm{max}}<2\alpha^2/(2-\alpha),
\end{equation}
and that $b(t, x)$ is bounded and  $\beta$-H\"older continuous
in $x$ uniformly in $t\in [0, 1]$,  for some
\begin{equation}\label{e:1.10}
\beta \in (\beta_0, 1]
\ \hbox{ with }\beta_0:=\alpha_{\mathrm{max}}/2+(\alpha_{\mathrm{max}}/\alpha
\1_{\{\alpha\leq 1\}} +\1_{\{\alpha>1\}})(1-\alpha).
\end{equation}
Then SDE \eqref{SDE1} has a unique strong solution for every $x\in \mR^d$,
where $d:=\sum_{j=1}^k d_j$, and \eqref{EQ9} holds.
\end{enumerate}
\end{corollary}

\medskip

Note that condition \eqref{e:1.9} implies that
$\alpha<2\alpha^2/(2-\alpha)$. The latter is equivalent to $\alpha> 2/3$.
If in Corollary \ref{C:1.3} (iii), $\alpha_j=\alpha$ for every $1\leq j\leq k$,
then  conditions \eqref{e:1.9} and \eqref{e:1.10} become
$$
\alpha>2/3  \quad \hbox{ and } \quad \beta \in (1-\alpha/2, 1],\ \hbox{ respectively.}
$$
An interesting open question is whether
the constraint $\alpha>2/3$ can be dropped.

For Corollary \ref{C:1.3} (iii),
let $\nu$ be the L\'evy measure of
the cylindrical stable process ${Z}$. We will in fact show in Example \ref{E:4.4} that,
when $\alpha=\min_{1\leq j\leq k}\alpha_j \in (1, 2)$,
condition {\bf \big(H$^{\alpha}_{\nu,K_0}$\big)} holds for some  $K_0>0$
but condition {\bf \big(H$^{\alpha^*}_{\nu,K_0}$\big)} fails for any $\alpha^*>\alpha$.
So  Hypothesis 1   of \cite{Pr1}
holds with this $\alpha$ for  the cylindrical stable process ${Z}$.
On the other hand, condition \eqref{EB1}
holds if and only if $2\gamma > \alpha_{\mathrm{max}}$.
Hence in the case $\alpha \in (1, 2)$,  Hypothesis 2 of \cite{Pr1}
fails when $\alpha_j$'s are not identical (i.e., when $\alpha_{\mathrm{max}}>\alpha$), and so the main results of \cite{Pr1}
are not applicable.

\medskip

The second main result of this paper is the following derivative formula of $\mE f(X_t(x))$.

\bt\label{Main2}
Under the assumptions of Theorem \ref{Main1}, if $Z_t=W_{S_t}$ is a
subordinate Brownian motion as described
in Example \ref{EEU} below,
then we have the following derivative formula:
\begin{align}
\nabla\mE f(X_t(x))=\mE\left[\frac{f(X_t(x))}{S_t}\int^t_0\nabla X_{s}(x)\dif
W_{S_s}\right],\ \ f\in C^1_b(\mR^d).\label{EQ15}
\end{align}
In particular, for any $p>1$, there is a constant $C_p>0$ such that for
any $f\in C^1_b(\mR^d)$ and $(t,x)\in(0,1)\times\mR^d$,
\begin{align}
|\nabla\mE f(X_t(x))|\leq C_p t^{-1/{\alpha}}(\mE |f(X_t(x))|^p)^{1/p}.
\label{EQ16}
\end{align}
\et

This paper is organized as follows: In Section 2, we solve a nonlocal
advection equation and obtain
estimates on the gradient of the solutions.
In particular, we derive a priori uniform $C^{1+\gamma}$ estimate
on the solution of the nonlocal advection equation.
Even when ${Z}$ is a rotationally symmetric stable process, our approach to the
a priori estimate is simpler and more elementary than that of \cite{Si2}.
In Section 3, we shall prove our main results by using Zvonkin's transform.
In Section 4, we give three examples to illustrate the main results of this paper,
from which Corollary \ref{C:1.3} follows.
In Appendix, we prove a continuous dependence result about the SDEs with jumps
with respect to the coefficients and the initial values.

\section{Differentiability of solutions of
nonlocal advection equations}

In this paper we use the following conventions.
The letter $C$ with or without subscripts
will denote a positive constant, whose value is not important
and may change from one appearance to another.
We write $f(x)\preceq g(x)$ to mean that there exists a constant
$C_0>0$ such that $f(x)\leq C_0 g(x)$; and $f(x)\asymp g(x)$ to mean
that there exist $C_1,C_2>0$ such that $C_1 g(x)\leq f(x)\leq C_2 g(x)$.

For a function $u(t, x)$ defined on $[0, 1]\times \mR^d$, sometimes we use $u_t (x)$ for $u(t, x)$.
Denote by $C^\infty_c (\mR^d)$ the space of smooth functions with compact
support on $\mR^d$.
For $\beta\in(0,1]$ and a function $f$ on $\mR^d$,
$$
[f]_\beta:=\sup_{x\not=y}\frac{|f(x)-f(y)|}{|x-y|^\beta},\quad
 \|f\|_\beta:= \|f\|_\infty+[f]_\beta,
$$
and for a function $f:[0,1]\times\mR^d\to\mR$,
$$
[f]_{\infty,\beta}:=\sup_{s\in[0,1]}[f_s]_\beta,  \quad
\|f\|_{\infty,\beta}:=\sup_{s\in[0,1]}\|f_s\|_\beta.
$$
Recall the following characterization for a H\"older continuous function $f$.
Let $P_\theta f$ be the Poisson integral of $f$ defined by
$$
P_\theta f(x):=\int_{\mR^d}f(y)p_\theta(x-y)\dif y,\quad \theta>0,
$$
where $p_\theta(x)$ is the density of a Cauchy process $Z_\theta$ given by
$$
p_\theta(x):=c_d\theta(\theta^2+|x|^2)^{-\frac{d+1}{2}}\asymp \theta(\theta+|x|)^{-d-1}.
$$
It is well-known (cf. \cite[Proposition 7 on p.142]{St}) that $\|f\|_\beta<\infty$
if and only if $f$ is bounded and
$$
\|\p_\theta P_\theta f\|_\infty\leq C\theta^{\beta-1} \quad \hbox{for every } \theta>0
$$
and
\begin{align}
\|f\|_\beta\asymp \|f\|_\infty+\sup_{\theta>0}\|\theta^{1-\beta}\p_\theta
P_\theta f\|_\infty.\label{EQ5}
\end{align}

The following commutator estimate result plays
an important role in our proof
of the H\"older regularity of the gradient in the case of $\alpha\in(0,1]$.
\bl\label{Le31}
For any $\beta,\gamma\in(0,1)$ with $\gamma\leq\beta$, there is a positive
constant $C=C(\beta,\gamma,d)$ such that
for any Borel functions $f, g$ on $\mR^d$,
$$
[\p_\theta P_\theta (fg)-f\p_\theta P_\theta g]_{\beta-\gamma}\leq C[f]_\beta
\, \|g\|_\infty \, \theta^{\gamma-1},\quad \theta>0,
$$
provided that $[f]_\beta$ and $\|g\|_\infty$ are finite. In particular,
if $g\equiv1$, then
$$
[\p_\theta P_\theta f]_{\beta-\gamma}\leq C[f]_\beta\theta^{\gamma-1},\ \ \theta>0.
$$
\el

\begin{proof}
It suffices to prove that
\begin{align}
&|\p_\theta P_\theta (fg)(x)-f\p_\theta P_\theta g(x)-\p_\theta P_\theta
(fg)(x')+f\p_\theta P_\theta g(x')|\no\\
&\leq C[f]_\beta\|g\|_\infty\theta^{\gamma-1}|x-x'|^{\beta-\gamma}.\label{EQ1}
\end{align}
By definition, we have
\begin{align}
\p_\theta P_\theta (fg)(x)-f\p_\theta P_\theta g(x)=\int_{\mR^d}(f(y)-f(x))g(y)
\p_\theta p_\theta(x-y)\dif y.\label{EQ6}
\end{align}
Notice the following easy estimates:
\begin{align}
|\p_\theta p_\theta(x)|\preceq (\theta+|x|)^{-d-1},\ \ |\nabla\p_\theta
p_\theta(x)|\preceq (\theta+|x|)^{-d-2},\label{EQ2}
\end{align}
and
\begin{align}
\int_{\mR^d}|x|^\beta(\theta+|x|)^{-d-k}\dif x\preceq \theta^{\beta-k},\quad  k\in\mN.\label{EQ3}
\end{align}
If $|x-x'|\geq \theta/2$, then (\ref{EQ1}) follows from
\begin{align*}
&\|\p_\theta P_\theta (fg)-f\p_\theta P_\theta g\|_\infty
\stackrel{(\ref{EQ6})}{\leq} [f]_\beta\|g\|_\infty\int_{\mR^d}|y|^\beta|
\p_\theta p_\theta(y)|\dif y\\
&\quad\stackrel{(\ref{EQ2})}{\preceq}[f]_\beta\|g\|_\infty
\int_{\mR^d}|y|^\beta(\theta+|y|)^{-d-1}\dif y\\
&\quad\stackrel{(\ref{EQ3})}{\preceq} [f]_\beta\|g\|_\infty
\theta^{\beta-1}\preceq [f]_\beta\|g\|_\infty\theta^{\gamma-1}|x-x'|^{\beta-\gamma}.
\end{align*}
Next, we assume
\begin{align}
|x-x'|\leq \theta/2.\label{EQ4}
\end{align}
Notice that
\begin{align*}
&\p_\theta P_\theta (fg)(x)-f\p_\theta P_\theta g(x)-(\p_\theta
P_\theta (fg)(x')-f\p_\theta P_\theta g(x'))\\
&\quad=\int_{\mR^d}(f(y)-f(x))g(y)(\p_\theta p_\theta(x-y)-
\p_\theta p_\theta(x'-y))\dif y\\
&\qquad+\int_{\mR^d}(f(x')-f(x))g(y)\p_\theta p_\theta(x'-y)\dif y=:I_1+I_2.
\end{align*}
For $I_1$, we have
\begin{align*}
|I_1|&\leq[f]_\beta\|g\|_\infty\int_{\mR^d}|x-y|^\beta|x-x'|
\left(\int^1_0|\nabla\p_\theta p_\theta(x-y+r(x'-x))|\dif r\right)\dif y\\
&\stackrel{(\ref{EQ2})}{\preceq}[f]_\beta\|g\|_\infty|x-x'|
\int_{\mR^d}|x-y|^\beta\left(\int^1_0(\theta+|x-y+r(x'-x)|)^{-d-2}\dif r\right)\dif y\\
&\stackrel{(\ref{EQ4})}{\preceq}[f]_\beta\|g\|_\infty|x-x'|
\int_{\mR^d}|x-y|^\beta(\theta+|x-y|)^{-d-2}\dif y\\
&\stackrel{(\ref{EQ3})}{\preceq}[f]_\beta\|g\|_\infty|x-x'|
\theta^{\beta-2}\stackrel{(\ref{EQ4})}{\preceq}[f]_\beta\|g
\|_\infty|x-x'|^{\beta-\gamma}\theta^{\gamma-1}.
\end{align*}
For $I_2$, we similarly have
\begin{align*}
|I_2|\preceq |x-x'|^\beta[f]_\beta \|g\|_\infty\int_{\mR^d}(\theta+|y|)^{-d-1}\dif y\preceq[f]_\beta\|g\|_\infty|x-x'|^{\beta-\gamma}\theta^{\gamma-1}.
\end{align*}
Combining the above estimates, we obtain (\ref{EQ1}).
\end{proof}
We also need the following lemma for treating the case of $\alpha\in(1,2)$.

\bl\label{Le22}
Suppose that {\bf \big(H$^{\alpha}_{\nu_1,K_0}$\big)} holds for
some $\alpha\in(1,2)$ and $K_0>0$. Then for any $\beta,\gamma\in[0,1]$,
there is a constant $K_1>0$ such that
$$
\|\nabla T^{\nu_1,0}_{t}f\|_\gamma\leq K_1
t^{(\beta-1-\gamma)/\alpha}
\|f\|_\beta
\quad \mbox{ for all } t\in(0,1).
$$
\el

\begin{proof}
Note that $\|\nabla T^{\nu_1,0}_{t}f\|_\infty\leq \|\nabla f\|_\infty$.
By \eqref{Es133} and the interpolation theorem,
we have
$$
\|\nabla T^{\nu_1,0}_{t}f\|_\infty\preceq t^{(\beta-1)/\alpha}\|f\|_\beta.
$$
On the other hand, by {\bf \big(H$^{\alpha}_{\nu_1,K_0}$\big)} we have
$$
\|\nabla^2T^{\nu_1,0}_{t}f\|_\infty=\|\nabla T^{\nu_1,0}_{t/2}\nabla T^{\nu_1,0}_{t/2}f\|_\infty
\preceq (t/2)^{-1/\alpha}\|\nabla T^{\nu_1,0}_{t/2}f\|_\infty \preceq (t/2)^{(\beta-2)/\alpha}\|f\|_\beta.
$$
Hence,
$$
[\nabla T^{\nu_1,0}_{t}f]_\gamma\leq 2\|\nabla^2 T^{\nu_1,0}_{t}f\|^\gamma_\infty\|\nabla
T^{\nu_1,0}_{t}f\|^{1-\gamma}_\infty\preceq t^{(\beta-1-\gamma)/\alpha}\|f\|_\beta.
$$
\end{proof}

For $\lambda\geq 0$, consider the following linear backward nonlocal parabolic system:
\begin{align}
\p_t u_t+(\sL_{\nu, \eta}-\lambda) u_t+b_t\cdot\nabla u_t+f_t=0,\quad \ u_1=0,\label{Eq1}
\end{align}
where $\sL_{\nu, \eta}$ is the infinitesimal generator of the L\'evy process ${Z}$,
and $b,f:[0,1]\times\mR^d\to\mR^d$ are bounded Borel
functions. Recalling decomposition \eqref{Dec}, we can write
\begin{align}
\sL_{\nu,\eta}=\sL_{\nu_0,0}+\sL_{\nu_1,0}+\sL_{\nu_2,\eta}.\label{Eq2}
\end{align}

The following theorem is the main result of this section.

\bt\label{Th31}
Suppose that  either {\bf \big(H$^{\alpha,\bar\alpha,\delta}_{\nu_1,K_0}$\big)}
holds for some $\alpha\in (0, 1]$,
$\bar\alpha,\delta\in (0, 1]$ and $K_0>0$ or
{\bf \big(H$^{\alpha}_{\nu_1,K_0}$\big)} holds for some $\alpha\in (1, 2)$
and $K_0>0$.
Suppose further that
$1-\alpha<\delta \bar \alpha $ in the case $\alpha\in (0, 1]$.
If for some $\beta \in ((1-\alpha)/\delta, \bar \alpha]$
in the case $\alpha \in (0, 1]$ and $\beta \in [0,1]$ in the case  $\alpha \in (1, 2)$,
it holds that
\begin{align}
\|b\|_{\infty,\beta}<\infty,\quad \|f\|_{\infty,\beta}<\infty,\label{EW22}
\end{align}
then for any
$\gamma\in (0,\beta-(1-\alpha)/\delta )$ when $\alpha\in(0,1]$ and
$\gamma \in (0, (\beta+\alpha-1)\wedge 1)$ when  $\alpha\in(1,2)$,
there exists a continuous function $u:[0,1]\times\mR^d\to\mR^d$ such that
for all $t\in [0, 1]$ and $\varphi\in C^\infty_c(\mR^d)$,
\begin{align}\label{UT}
\<u_t,\varphi\>=\int^t_0\<u_s,(\sL^*_{\nu, \eta}-\lambda)\varphi\>\dif s+
\int^t_0\<b_s\cdot\nabla u_s,\varphi\>\dif s+\int^t_0\<f_s,\varphi\>\dif s
\end{align}
with
\begin{align}\label{Es101}
\sup_{t\in[0, 1]}\|u_t(\cdot)\|_\infty\leq \sup_{t\in[0, 1]}\|f_t(\cdot)\|_\infty,
\end{align}
and for some $\theta_0>0$ and all $\lambda\geq 0$,
\begin{align}\label{Es01}
\|\nabla u\|_{\infty,\gamma}\leq C(1\vee\lambda)^{-\theta_0}\|f\|_{\infty,\beta}.
\end{align}
Here
$C=C(d, \alpha,\beta, K_0,   \bar \alpha, \delta, \|b\|_{\infty,\beta},
\gamma,|\nu_0|(\mR^d))$,
$\<u,\varphi\>:=\int u\varphi\dif x$ and $\sL^*_{\nu,\eta}$ is the adjoint operator of $\sL_{\nu,\eta}$.
\et

We will first prove several lemmas before we present the proof of the theorem above.
Let $Z^{(1)}$ and $Z^{(2)}$
be two independent L\'evy processes
with generators $\sL_{\nu_0+\nu_1,0}$
and $\sL_{\nu_2,\eta}$. Clearly,
\begin{align}
Z_t\stackrel{(d)}{=}Z^{(1)}_t+Z^{(2)}_t.\label{Es98}
\end{align}

\bl\label{Le24}
Assume that $b,f\in L^\infty([0,1];C^\infty_b(\mR^d;\mR^d))$.
There exists a unique solution
$u_t(x)\in C([0,1];C^\infty_b(\mR^d;\mR^d))$
to equation \eqref{Eq1} with the following probabilistic representation:
\begin{align}
u_t(x)=\int^1_t\e^{\lambda(t-s)}\mE f_s(X_{t,s}(x))\dif s,\label{ET1}
\end{align}
where $X_{t,s}(x)=X_{t,s}$ is the unique solution to the following SDE:
\begin{align}
X_{t,s}=x+\int^s_t b_r(X_{t,r})\dif r+Z_s-Z_t,\quad s\geq t.\label{ET2}
\end{align}
Moreover, we have the following a priori estimate
\begin{align}
\sup_{t\in [0, 1]}\|u_t\|_\infty\leq \sup_{t\in [0, 1]}\|f_t\|_\infty.\label{Ma0}
\end{align}
\el

\begin{proof}
The existence and uniqueness of $u_t(x)$ and the representation (\ref{ET1}) follow from
\cite[Theorem 4.4]{Zh2}. The estimate
(\ref{Ma0}) immediately follows from (\ref{ET1}).
\end{proof}

\bl\label{Le25}
Suppose that {\bf \big(H$^{\alpha,\bar\alpha,\delta}_{\nu_1,K_0}$\big)}
 holds for some $\alpha\in (0, 1]$, $\bar\alpha,\delta\in (0, 1]$
 and $K_0>0$, and that $b,  f\in L^\infty([0,1];C^\infty_b(\mR^d;\mR^d))$.
Suppose further that
$1-\alpha<\delta\bar\alpha$ in the case $\alpha\in (0, 1]$.
Let $u$ be the solution of \eqref{Eq1}.
Then for any
$\beta_1,\beta_2\in((1-\alpha)/\delta,\bar\alpha]$,
there is a constant $C>0$   depending only on
$K_0,\alpha,\delta,\beta_1,\beta_2$, $[b]_{\infty,\beta_1}$  and
$|\nu_0|(\mR^d)$ such that for all
 $\lambda\geq 0$,
\begin{align}
\sup_{t\in[0, 1]}\|\nabla u_t\|_\infty\leq C(1\vee\lambda)^{(1-\alpha- \delta\beta_2)/{\alpha}}
[f]_{\infty,\beta_2}.\label{Es001}
\end{align}
\el

\begin{proof}
(i) We first assume that $\eta=0$ and that $\nu_2=0$ in decomposition (\ref{Dec}).
Fix $x_0\in\mR^d$ and let $y_t$ satisfy the following ODE:
$$
\dot y_t=-b_t(x_0+y_t) \ \hbox{ with }  y_0=0.
$$
Define
\begin{align}
\tilde u_t(x):=u_t(x+x_0+y_t),\quad  \tilde f_t(x):=f_t(x+x_0+y_t)\label{def1}
\end{align}
and
$$
\tilde b_t(x):=b_t(x+x_0+y_t)-b_t(x_0+y_t).
$$
Clearly, by (\ref{Eq1}) and (\ref{Eq2}), $\tilde u$ satisfies
\begin{align*}
\p_t \tilde u_t+(\sL_{\nu_1,0}-\lambda) \tilde u_t+\tilde b_t\cdot\nabla
\tilde u_t+\sL_{\nu_0,0}\tilde u_t +\tilde f_t =0,\quad \tilde u_1=0.
\end{align*}
We have by the representation \eqref{ET1} (with $b=0$ there)
$$
\tilde u_t(x)=\int^1_t\e^{\lambda(t-s)}{T}^{\nu_1,0}_{s-t}
\left(  \tilde b_s\cdot \nabla \tilde u_s +
\sL_{\nu_0,0}\tilde u_s +  \tilde f_s\right)(x)\dif s.
$$
Fix $\beta_1,\beta_2\in ((1-\alpha)/\delta, \bar\alpha]$.
Note that
by the definition of $\tilde b_s$,
$$
|\tilde b_s(y)\cdot \nabla \tilde u_s(y)|\leq \|\nabla \tilde
u_s(\cdot)\|_\infty[b_s(\cdot)]_{\beta_1}|y|^{\beta_1}
\ \ \hbox{for all } y\in\mR^d ,
$$
and that
$$
| \nabla T^{\nu_1,0}_{s-t} (\sL_{\nu_0,0}\tilde u_s )(x)|
\leq \| \nabla  (\sL_{\nu_0,0}\tilde u_s )\|_\infty
\leq 2 |\nu_0|(\mR^d) \| \nabla \tilde u_s(\cdot) \|_\infty.
$$
We have by \eqref{Es13} that for $t\in [0, 1]$,
\begin{align*}
|\nabla \tilde u_t(0)|&\leq \int^1_t\e^{\lambda(t-s)}
\Big(K_0[b_s(\cdot)]_{\beta_1}(s-t)^{(\delta\beta_1-1)/{\alpha}}+
2|\nu_0|(\mR^d)\Big)\|\nabla \tilde u_s(\cdot)\|_\infty\dif s\\
&\quad+K_0\int^1_t\e^{\lambda(t-s)}(s-t)^{ (\delta\beta_2-1)/{\alpha}}[\tilde f_s(\cdot)]_{\beta_2}\dif s.
\end{align*}
By \eqref{def1} and the arbitrariness of $x_0$, one in fact has
\begin{align}
\|\nabla u_t(\cdot)\|_\infty\leq C\int^1_t(s-t)^{(\delta\beta_1-1)/{\alpha}}
\|\nabla u_s(\cdot)\|_\infty\dif s
 +C[f]_{\infty,\beta_2}(1\vee\lambda)^{(1-\alpha-\delta\beta_2)/{\alpha}}.\label{In9}
\end{align}
By Gronwall's inequality, we obtain (\ref{Es001}).

(ii) Next we consider the general case. Fix
$t_0\in[0,1)$ and a c\'adl\'ag function $\ell:[0,1]\to\mR^d$, and define
$$
b^\ell_r(x):=b_r(x-\ell_{t_0}+\ell_r),\ \ f^\ell_r(x):=f_r(x-\ell_{t_0}+\ell_r).
$$
Let $Y^\ell_{t,s}(x):=Y^\ell_{t,s}$ be the solution to the following SDE:
$$
Y^\ell_{t,s}=x+\int^s_t b^\ell_r(Y^\ell_{t,r})\dif r+Z^{(1)}_s-Z^{(1)}_t,\ s\geq t.
$$
Since $Z^{(1)}$ and $Z^{(2)}$ are independent,
by \eqref{Es98} and
the uniqueness in law of the solution to SDE \eqref{ET2}, we have
$$
X_{t_0,\cdot}(x)\stackrel{(d)}{=}Y^{Z^{(2)}}_{t_0,\cdot}(x)-Z^{(2)}_{t_0}+Z^{(2)}_\cdot,
$$
and so by  (\ref{ET1}),
$$
u_{t_0}(x)=\mE\left(\int^1_{t_0}\e^{\lambda(t_0-s)}\mE\Big[ f^\ell_s
(Y^\ell_{t_0,s}(x))\Big] \dif s\Big|_{\ell=Z^{(2)}}\right).
$$
Now we define
$$
u^\ell_t(x):=\int^1_t\e^{\lambda(t-s)}\mE\Big[
f^\ell_s(Y^\ell_{t,s}(x))\Big] \dif s.
$$
Then by Lemma \ref{Le24}, $u^\ell_t(x)$ is a solution to the following equation:
$$
\p_t u^\ell+(\sL_{\nu_0+\nu_1,0}-\lambda) u^\ell+b^\ell\cdot\nabla u^\ell
+f^\ell=0,\ \ u^\ell_1=0.
$$
In view of
$$
[b^\ell]_{\infty,\beta_1}=[b]_{\infty,\beta_1},\ \
[f^\ell]_{\infty,\beta_2}=[f]_{\infty,\beta_2},
$$
by what has been proved in (i),
we have for any c\'adl\'ag function $\ell$,
$$
\|\nabla u^\ell\|_\infty\leq
C_2(1\vee\lambda)^{(1-\alpha-\delta\beta_2)/{\alpha}}[f]_{\infty,\beta_2},
$$
which in turn gives (\ref{Es001}) by noting that $\nabla u_{t_0}(x)=
\mE\Big[ \nabla u^\ell_{t_0}(x)|_{\ell=Z^{(2)}}\Big]$ and $t_0$ is arbitrary.
\end{proof}

\bl
Suppose that {\bf \big(H$^{\alpha,\bar\alpha,\delta}_{\nu_1,K_0}$\big)}
 holds for some $\alpha\in (0, 1]$, $\bar\alpha,\delta\in (0, 1]$
 and $K_0>0$, and that $b,  f\in L^\infty([0,1];C^\infty_b(\mR^d;\mR^d))$.
Suppose further that
$1-\alpha<\delta\bar\alpha$ in the case $\alpha\in (0, 1]$.
Let $u$ be the solution of \eqref{Eq1}.
Then for any
$\beta\in((1-\alpha)/\delta,\bar\alpha]$
and $\gamma\in(0, \beta - (1-\alpha)/\delta)$,
there exists a constant $C>0$ depending only on
$d,\alpha, \bar \alpha, \delta, K_0, \beta, \gamma $,
$[b]_{\infty,\beta}$  and $|\nu_0|(\mR^d)$ such
that for all $\lambda\geq 0$,
\begin{align}
[\nabla u]_{\infty,\gamma}\leq C(1\vee\lambda)^{(1-\alpha-\delta(\beta-\gamma))/{\alpha}}
[f]_{\infty,\beta}.\label{Es011}
\end{align}
\el

\begin{proof}
Fix $\gamma\in(0, \beta - (1-\alpha)/\delta)$. For $\theta>0$, define
$$
w^\theta_t(x):=\p_\theta P_\theta u_t(x)
$$
and
$$
g^\theta_t(x):=\p_\theta P_\theta(b_t\cdot\nabla u_t)(x)-b_t(x)
\cdot\nabla \p_\theta P_\theta u_t(x)+\p_\theta P_\theta f_t(x),
$$
then
$$
\p_t w^\theta_t+(\sL_{\nu, \eta}-\lambda)w^\theta_t+b_t\cdot\nabla w^\theta_t+g^\theta_t=0,
\quad \ w^\theta_1=0.
$$
Since $\beta-\gamma> (1-\alpha)/\delta$, by (\ref{Es001}) with $\beta_1=\beta$ and $\beta_2=\beta-\gamma$,
we have
$$
\sup_{t\in [0, 1]}\|\nabla w^\theta_t(\cdot)\|_\infty\leq
C(1\vee\lambda)^{ (1-\alpha-\delta(\beta-\gamma))/{\alpha}}[g^\theta]_{\infty,\beta-\gamma},
$$
and by Lemma \ref{Le31},
\begin{align*}
[g^\theta]_{\infty,\beta-\gamma}&\preceq[b]_{\infty,\beta}
\sup_{t\in [0, 1]}\|\nabla u_t\|_{\infty}
\theta^{\gamma-1}+[f]_{\infty,\beta}\theta^{\gamma-1}
\stackrel{(\ref{Es001})}{\preceq} [f]_{\infty,\beta}\theta^{\gamma-1}.
\end{align*}
Hence,
$$
\sup_{t\in [0, 1]}\|\p_\theta P_\theta \nabla u_t\|_\infty \leq C
(1\vee\lambda)^{ (1-\alpha-\delta(\beta-\gamma))/{\alpha}}
[f]_{\infty,\beta}\theta^{\gamma-1}
\mbox{ for every $\theta>0$,}
$$
which yields \eqref{Es011} by (\ref{EQ5}) and (\ref{Es001}).
\end{proof}

\bl
Suppose that {\bf \big(H$^{\alpha}_{\nu_1,K_0}$\big)}
 holds for some $\alpha\in (1, 2)$
 and $K_0>0$, and that $b,  f\in L^\infty([0,1];C^\infty_b(\mR^d;\mR^d))$.
Let $u$ be the solution of \eqref{Eq1}.
Then for any $\beta\in[0, 1]$ and $\gamma\in[0, (\beta+\alpha-1)\wedge 1)$,
there exists a constant $C>0$ depending only on $d, \alpha, K_0,
\gamma,\beta$, $\|b\|_{\infty,\beta}$  and $|\nu_0|(\mR^d)$ such
that for all $\lambda\geq 0$,
\begin{align}
\|\nabla u\|_{\infty,\gamma}\leq C(1\vee\lambda)^{(1-\alpha-\beta+\gamma)/{\alpha}}
\|f\|_{\infty,\beta}.\label{Es0111}
\end{align}
\el

\begin{proof}
As in the proof of Lemma \ref{Le25}, we first assume that $\eta=0$ and that $\nu_2=0$
in decomposition (\ref{Dec}).
By the representation \eqref{ET1} (with $b=0$ there), we have
$$
u_t(x)=\int^1_t\e^{\lambda(t-s)}{T}^{\nu_1,0}_{s-t}\left(b_s\cdot \nabla u_s +\sL_{\nu_0,0}u_s
+ f_s\right)(x)\dif s.
$$
Without loss of generality, we assume
$\gamma\in[\beta, (\beta+\alpha-1)\wedge 1)$.
By Lemma \ref{Le22}, we have
\begin{align*}
\|\nabla u_t\|_\gamma&\preceq\int^1_t\e^{\lambda(t-s)}(s-t)^{\frac{\beta-1-\gamma}{\alpha}}
\left(\|b_s\cdot \nabla u_s\|_\beta +\|f_s\|_\beta\right)\dif s
+\int^1_t\e^{\lambda(t-s)}\|\nabla{T}^{\nu_1,0}_{s-t}\sL_{\nu_0,0}u_s\|_\gamma\dif s\\
&\preceq\int^1_t\e^{\lambda(t-s)}(s-t)^{(\beta-1-\gamma)/{\alpha}}\left(\|b_s\|_\beta\|\nabla u_s\|_\beta +\|f_s\|_\beta\right)\dif s
+\int^1_t\e^{\lambda(t-s)}\|\nabla u_s\|_\gamma\dif s\\
&\preceq\int^1_t((s-t)^{(\beta-1-\gamma)/{\alpha}}+1)\|\nabla u_s\|_\gamma\dif s
+\|f\|_{\infty,\beta}\int^1_t\e^{\lambda(t-s)}(s-t)^{(\beta-1-\gamma)/{\alpha}}\dif s,
\end{align*}
which yields
\eqref{Es0111} by Gronwall's inequality.
For the general case, we can follow the same argument in (ii) of Lemma \ref{Le25} to derive \eqref{Es0111}.
\end{proof}

Now we are ready to give

\begin{proof}[Proof of Theorem \ref{Th31}]
Suppose that $b$ and $f$ satisfy \eqref{EW22}.
Let $\varrho$ be a non-negative smooth function with compact support in $\mR^d$
having $\int_{\mR^d} \varrho (x) \dif x=1$. For $n\in\mN$, define $\varrho_n (x):=n^d \varrho (nx)$ and
\begin{align}
b^n_t:=\varrho_n* b_t,\ \ f^n_t:=\varrho_n*f_t.\label{EQ10}
\end{align}
Clearly, $b^n, f^n\in L^\infty([0,1];C^\infty_b(\mR^d, \mR^d))$ and
$$
\|b^n\|_{\infty,\beta}\leq \|b\|_{\infty,\beta},\ \ \|f^n\|_{\infty,\beta}\leq \|f\|_{\infty,\beta}.
$$
Let $u^n_t$ be the solution to the following equation:
\begin{align}
\p_tu^n_t+(\sL_{\nu,\eta}-\lambda)u^n_t+b^n_t\cdot\nabla u^n_t+f^n_t=0, \quad u^n_1=0.\label{EA1}
\end{align}
By \eqref{Es001}, (\ref{Es011}) and  \eqref{Es0111}, there is an $\theta_0>0$ such that for all $\lambda\geq 0$,
\begin{align}
\sup_{n}\|\nabla u^n\|_{\infty,\gamma}\leq C(1\vee\lambda)^{-\theta_0}\|f\|_{\infty,\beta},\label{ET3}
\end{align}
and by (\ref{Ma0}),
\begin{align}
\sup_{t\in[0,1]}\|u^n_t\|_\infty\leq \sup_{t\in[0,1]}\|f^n_t\|_\infty\leq \sup_{t\in[0,1]}\|f_t\|_\infty.
\label{ET4}
\end{align}
Moreover, by the representation \eqref{ET1} (with $b=0$ there),
we can write
$$
u^n_t(x)=\int^1_t\e^{\lambda(t-s)}{T}^{\nu,\eta}_{s-t}\left(b^n_s\cdot\nabla u^n_s+f_s^n\right)(x)\dif s.
$$
Using (\ref{ET3}) and (\ref{ET4}), one can easily show that for any $R>0$,
$$
\lim_{|t-t'|\to 0}\sup_{n}\sup_{|x|\leq R}|u^n_t(x)-u^n_{t'}(x)|=0.
$$
Hence, by Ascoli-Arzela's lemma, there is a subsequence (still
denoted by $u^
n$) and a function $u$ with
$$
\|\nabla u\|_{\infty,\gamma}\leq C(1\vee\lambda)^{-\theta_0}\|f\|_{\infty,\beta},
\quad \sup_{t\in [0, 1]}\|u_t\|_\infty\leq \sup_{t\in [0, 1]}\|f_t\|_\infty
$$
such that
\begin{align}
\lim_{n\to\infty}\sup_{t\in[0,1]}\sup_{|x|\leq R}|u^n_t(x)-u_t(x)|=0, \quad \forall R>0.\label{EQ11}
\end{align}
On the other hand, noticing the following interpolation inequality (cf. \cite[Theorem 3.2.1]{Kr})
$$
\|\nabla \phi \|_\infty\leq C\|\nabla \phi \|^{{1}/{1+\gamma}}_\gamma\|\phi\|^{{\gamma}/{1+\gamma}}_\infty,
$$
by \eqref{ET3} and \eqref{EQ11}, we further have
\begin{align}
\lim_{n\to\infty}\sup_{t\in[0,1]}\sup_{|x|\leq R}|\nabla u^n_t(x)-\nabla u_t(x)|=0
\quad \hbox{for every } R>0.\label{EQ111}
\end{align}
Thus  by \eqref{EA1}, \eqref{EQ11} and
\eqref{EQ111},
it is easy to see that $u$ satisfies (\ref{UT}).
\end{proof}

\bc\label{Re1}
Under the assumption of Theorem \ref{Th31},
if we further assume that for some
$\gamma_0\in(0, \beta-(1-\alpha)/\delta)$ in the case of $\alpha\in(0,1]$ and
$\gamma_0\in(0, (\beta+\alpha-1)\wedge 1)$
in the case of $\alpha\in(1,2)$,
\begin{align}
\int_{|z|\leq 1}|z|^{1+\gamma_0}\nu(\dif z)<\infty, \label{EB2}
\end{align}
then the solution $u$ of equation \eqref{Eq1} satisfying  \eqref{Es101}  and \eqref{Es01}
for some $\gamma >\gamma_0$
is  a classical
solution; that is, $\sL_{\nu, \eta} u_s (x)$ and $\nabla u_s (x)$
exists pointwise and is continuous in $x$, and for
all $t\in[0,1]$ and $x\in\mR^d$,
\begin{align}
u_t(x)+\int^1_t(\sL_{\nu,\eta}-\lambda)u_s(x)\dif s+\int^1_tb_s(x)
\cdot\nabla u_s(x)\dif s+\int^1_tf_s(x)\dif s=0.\label{Eq9}
\end{align}
\ec

\begin{proof}
Since
$\|\nabla u\|_{\infty,\gamma}<\infty$
for some $\gamma\in(\gamma_0,\beta-(1-\alpha)/\delta)$ in the case of $\alpha\in(0,1]$ and
$\gamma\in(\gamma_0,(\beta+\alpha-1)\wedge 1)$ in the case of $\alpha\in(1,2)$,
by \eqref{EB2},
it is easy to check that
$$
x\mapsto \sL_{\nu,\eta} u_s(x)\mbox{ is continuous.}
$$
Hence, by (\ref{UT}), equation (\ref{Eq9}) is satisfied for
all $t\in[0,1]$ and $x\in\mR^d$.
\end{proof}

\section{Stochastic flow and Bismut formula}

Suppose that either {\bf \big(H$^{\alpha,\bar\alpha,\delta}_{\nu_1,K_0}$\big)} holds for some $\alpha\in (0, 1]$,
$\bar\alpha,\delta\in (0, 1]$ and $K_0>0$ or
{\bf \big(H$^{\alpha}_{\nu_1,K_0}$\big)} holds for some $\alpha\in (1, 2)$
and $K_0>0$. Suppose also that
\eqref{EB1} and \eqref{Ew22}  hold for some
$$
\gamma \in (0, 1) \ \hbox{with} \ \gamma +(1-\alpha)/\delta<\bar\alpha \ \hbox{and} \ \beta \in (\gamma + (1-\alpha)/\delta, \bar\alpha]
\mbox{ in the case of } \alpha\in(0,1]
$$
and
$$
\gamma \in (0, 1) \ \hbox{and} \ \beta
\in ((\gamma + 1-\alpha)^+, 1]
\hbox{ in the case of } \alpha\in(1,2).
$$
Notice that (\ref{EB1}) implies (\ref{EB2}) with $\gamma_0=\gamma$.
Hence, for $\lambda\geq 0$, by Corollary \ref{Re1}, the following
nonlocal equation has a classical solution $u$:
\begin{align*}
\p_t u_t+(\sL_{\nu,\eta}-\lambda) u_t+b_t\cdot\nabla u_t+b_t=0,\ \ u_1(x)=0.
\end{align*}
Similarly, let $b^n$ be defined by (\ref{EQ10}) and let $u^n$ be the solution to the following equation:
$$
\p_t u^n_t+(\sL_{\nu,\eta}-\lambda) u^n_t+b^n_t\cdot\nabla u^n_t+b^n_t=0,
\quad u^n_1(x)=0.
$$
Using the same argument leading to
\eqref{EQ11} and \eqref{EQ111},
we see that there is a subsequence, still denoted by $u^n$, such that
\begin{align}
\lim_{n\to\infty}\sup_{t\in[0,1]}\sup_{|x|\leq R}|\nabla^j u^n_t(x)-\nabla^ju_t(x)|=0
\quad \hbox{for every } R>0
\hbox{ and } j=0,1.\label{Lim}
\end{align}
For simplicity, we use the following convention:
$$
u^\infty:=u,\quad \ b^\infty:=b,\quad \ \mN_\infty:=\mN\cup\{\infty\}.
$$
By (\ref{Es01}), one can choose $\lambda$ sufficiently large, independent of $n\in\mN_\infty$, such that
\begin{align}\label{Lp8}
\|\nabla u^n_t(\cdot)\|_{\infty}+\sup_{x\not=x'}\frac{|\nabla u^n_t(x)-
\nabla u^n_t(x')|}{|x-x'|^\gamma}\leq\frac{1}{2}
\quad \hbox{for every } n\in\mN_\infty \hbox{ and } t\in [0, 1].
\end{align}
Define
$$
\Phi^n_t(x)=x+u^n_t(x),  \quad n\in\mN_\infty.
$$
Since for each $t\in[0,1]$,
$$
\frac{1}{2}|x-y|\leq|\Phi^n_t(x)-\Phi^n_t(y)|\leq\frac{3}{2}|x-y|,
$$
$x\mapsto\Phi^n_t(x)$ is a diffeomorphism
with
\begin{align} \label{Op1}
1/2\leq |\nabla\Phi^n_t(x)|\leq {3}/{2}
\quad \hbox{and} \quad |\nabla (\Phi^{n}_t)^{-1} (x)|\leq 2,
\end{align}
where
$(\Phi^{n}_t)^{-1}$
denotes the inverse function of $x\mapsto \Phi^n_t(x)$.

\bl
There is a constant $C=C(d)>0$ such that for all $t\in[0,1]$ and $n\in\mN_\infty$,
\begin{align}\label{EQ7}
\|\nabla\Phi^n_t\|_\gamma+\|\nabla (\Phi^{n}_t)^{-1}\|_\gamma\leq C.
\end{align}
Moreover, for each $t\in[0,1]$, $R>0$ and $j=0,1$, we have
\begin{align}
\lim_{n\to\infty}\sup_{t\in[0,1]}\sup_{|x|\leq R} \left| \nabla^j\Phi^n_t(x)-\nabla^j
\Phi^\infty_t(x) \right|=0,\label{Lim1}
\end{align}
and
\begin{align}\label{Lim2}
\lim_{n\to\infty}\sup_{t\in[0,1]}\sup_{|x|\leq R}
\left| \nabla^j (\Phi^{n}_t)^{-1} (x)-\nabla^j (\Phi^{\infty}_t)^{-1}(x) \right|=0.
\end{align}
\el

\begin{proof}
(i)
For notational simplicity, we   drop the superscript ``$n$''.
Clearly,
$\sup_{t\in [0, 1]}\|\nabla\Phi_t(\cdot)\|_\gamma < d+1$.
In view of
$$
(\nabla\Phi_s)^{-1}(x)-(\nabla\Phi_s)^{-1}(x')
=(\nabla\Phi_s)^{-1}(x) \left(\nabla\Phi_s(x')-\nabla\Phi_s(x) \right)(\nabla\Phi_s)^{-1}(x'),
$$
we have by \eqref{Lp8} and \eqref{Op1},
$$
[(\nabla\Phi_s)^{-1}]_\gamma\leq \|(\nabla\Phi_s)^{-1}\|_\infty^2
[\nabla\Phi_s]_\gamma
=\|(\nabla\Phi_s)^{-1}\|_\infty^2
 [\nabla u_s]_\gamma \leq 2
 \quad \mbox{ for all } s\in [0, 1].
$$
Hence by \eqref{Op1} again, for all $s\in [0, 1]$,
$$
\|\nabla\Phi^{-1}_s\|_\gamma
=\|(\nabla\Phi_s)^{-1}(\Phi^{-1}_s)\|_\gamma\leq
\|(\nabla\Phi_s)^{-1}\|_\infty+\|\nabla\Phi_s^{-1}\|^\gamma_\infty
[(\nabla\Phi_s)^{-1}]_\gamma\leq
2+2^{\gamma +1}.
$$

(ii) Properties (\ref{Lim1}) and (\ref{Lim2}) follow from the definitions of $\Phi_t$ and $\Phi^{-1}_t$,
\eqref{Lim} and \eqref{Lp8}.
\end{proof}

\medskip

For any given $n\in\mN_\infty$, define
\begin{align}\label{Lp00}
g^n_s(y,z):=\Phi^n_s \left( (\Phi^{n}_s)^{-1} (y)+z \right)-y.
\end{align}

\bl\label{L:3.2}
For $\gamma_1,\gamma_2\geq0$ with $\gamma_1+\gamma_2 =\gamma \in (0, 1)$,
there is a positive constant
 $C_1 =C_1 (d, \lambda, \gamma_1,\gamma_2)$ such that
for all $n\in\mN_\infty$, $t\in[0,1]$ and $y, z\in\mR^d$,
\begin{equation}\label{e:3.8}
\|\nabla g^n(\cdot,z)\|_{\infty,\gamma_1}
\leq C_1 (1\wedge|z|^{\gamma_2})
\quad \hbox{and} \quad  |g^n_s(y,z)|\leq 3|z|/2.
\end{equation}
Moreover, for each $t\in[0,1]$, $y, z\in\mR^d$ and $j=0,1$, we have
\begin{equation}\label{e:3.9}
\lim_{n\to\infty}\nabla^j_yg^n_t(y,z)=\nabla^j_yg^\infty_t(y,z).
\end{equation}
\el

\begin{proof}
For notational simplicity, we drop the superscript ``$n$''.
Since
$$
\nabla_y g_s(y,z)=\nabla\Phi_s \left(\Phi^{-1}_s(y)+z \right) \cdot\nabla\Phi^{-1}_s(y)-\mI,
$$
we have
$$
\|\nabla g_s(\cdot,z)\|_\infty\leq
2 \| \nabla\Phi_s \|_\gamma
(1\wedge |z|^\gamma)
\|\nabla\Phi^{-1}_s\|_\infty\stackrel{(\ref{EQ7})}{\leq}
C(1\wedge |z|^\gamma)
$$
and
\begin{align*}
[\nabla g_s(\cdot,z)]_\gamma&\leq [\nabla\Phi_s(\Phi^{-1}_s
(\cdot)+z)]_\gamma\|\nabla\Phi^{-1}_s\|_\infty+\|\nabla\Phi_s\|_\infty[\nabla\Phi^{-1}_s]_\gamma\\
&\leq [\nabla\Phi_s]_{\gamma}\|\nabla\Phi^{-1}_s\|^{1+\gamma}_\infty
+\|\nabla\Phi_s\|_\infty[\nabla\Phi^{-1}_s]_{\gamma}\stackrel{(\ref{EQ7})}{\leq} C.
\end{align*}
Thus, by definition, for $\gamma_1+\gamma_2=\gamma$, we have
\begin{align*}
[\nabla g_s(\cdot,z)]_{\gamma_1}\leq (2\|\nabla g_s(\cdot,z)\|_\infty)^{\gamma_2/\gamma}
[\nabla g_s(\cdot,z)]^{\gamma_1/\gamma}_\gamma\leq C(1\wedge|z|^{\gamma_2}),
\end{align*}
which in turn gives the first estimate in \eqref{e:3.8}.
The second inequality in  \eqref{e:3.8} follows from \eqref{Op1}
and the definition  of  $g^n$ .
Property \eqref{e:3.9} follows from \eqref{Lim1}, \eqref{Lim2}, and the definition  of  $g^n$.
\end{proof}

Taking $\gamma_1=0$ in Lemma \ref{L:3.2} yields that there is a constant
$C_0=C_0(d, \lambda,  \gamma)>0$ so that
\begin{equation}\label{e:3.10}
\|\nabla g^n(\cdot,z)\|_{\infty}
\leq C_0 (1\wedge|z|^{\gamma })
\quad \hbox{and} \quad  |g^n_s(y,z)|\leq 3|z|/2.
\end{equation}
Choose $r_0\in (0, 1)$ so that
\begin{equation}\label{e:3.11}
C_0 r_0^\gamma + 3r_0/2<1.
\end{equation}
Such a choice of $r_0$ will be used below to establish the
$C^1$-stochastic diffeomorphic property of
the unique solution $Y^n$ of SDE \eqref{Eq11}.

Let $N(\dif t,\dif z)$ be the Poisson random measure associated with ${Z}$, i.e.,
$$
N((0,t]\times\Gamma):=\sum_{0<s\leq t}1_\Gamma(Z_s-Z_{s-}), \quad t>0, \
\Gamma\in {\mathcal B}(\mR^d\setminus\{0\}).
$$
Let $\tilde N(\dif t,\dif z):=N(\dif t,\dif z)-\dif t\nu(\dif z)$ be the
compensated Poisson random measure.
By the L\'evy-It\^o decomposition, we can write for
each $r>0$,
$$
Z_t=\int^t_0\!\!\!\int_{|z|<r}z\tilde N(\dif s,\dif z)+\int^t_0\!\!\!
\int_{|z|\geq r}z N(\dif s,\dif z) +
\eta_r t,
$$
where $\eta_r \in \mR^d$ is a constant vector depending on $r$.

Recall that $r_0\in (0, 1)$ is the constant in \eqref{e:3.11}.
For any given $n\in\mN_\infty$, define
\begin{align}
 a^n_s(y):= \eta_{r_0} +
\lambda u^n_s \big( (\Phi^{n}_s)^{-1}(y)\big) -\int_{|z|\geq {r_0}}
\left( u^n_s \big( (\Phi^{n}_s)^{-1} (y)+z \big)-
u^n_s ( (\Phi^{n}_s)^{-1}(y))\right) \nu(\dif z) .
\label{Lp0}
\end{align}
We have

\bl\label{Le33}
There is a  positive constant
$C_2 =C_2 (d,  \lambda , \gamma, r_0)$   such that
for all $n\in\mN_\infty$, $t\in[0,1]$ and
$y\in\mR^d$,
\begin{align}
\|\nabla a^n\|_{\infty,\gamma}\leq C_2 \quad
\hbox{and} \quad  |a^n_s(y)|\leq C_2 (1+\| b\|_\infty)  . \label{EK1}
\end{align}
Moreover, for each $t\in[0,1]$,
$y\in\mR^d$
and $j=0,1$, we have
\begin{align}
\lim_{n\to\infty}\nabla^j_y a^n_t(y)=\nabla^j_ya^\infty_t(y).\label{Lim4}
\end{align}
\el

\begin{proof}
 For notational simplicity, we drop the superscript ``$n$''.
Since
$$
\nabla( u_s(\Phi^{-1}_s))=(\nabla u_s)(\Phi^{-1}_s)\cdot\nabla\Phi^{-1}_s,
$$
we have by (\ref{EQ7}) that for all $s\in [0, 1]$,
\begin{align*}
\|\nabla( u_s(\Phi^{-1}_s))\|_\gamma&\leq \|\nabla u_s(\Phi^{-1}_s)
\|_\gamma\|\nabla\Phi_s^{-1}\|_\infty+\|\nabla u_s\|_\infty
\|\nabla\Phi_s^{-1}\|_\gamma\\
&\leq \|\nabla u_s\|_\gamma\|\nabla
\Phi^{-1}_s\|_\infty (1+ \|\Phi^{-1}_s\|^{\gamma}_\infty)
+\|\nabla u_s\|_\infty\|\nabla\Phi_s^{-1}\|_\gamma\leq C.
\end{align*}
Hence $\|\nabla a\|_{\infty,\gamma}\leq C_1$ by \eqref{Lp0}.
The second inequality in \eqref{EK1}
follows from the definition of $a_s(y)$ and the fact that $u^n$ is uniformly bounded due to \eqref{Es101} of Theorem \ref{Th31}.
  Property \eqref{Lim4} follows from
 \eqref{Lim}, \eqref{Lp8}, \eqref{Lim1}, \eqref{Lim2},
 and the definition of $a^n$.
\end{proof}

The following lemma is a direct application of It\^o's formula.
\bl\label{Le3}
Let $\Phi^n_t(x)$ be defined as above. For $n\in\mN_\infty$, $X^n_t$
satisfies
\begin{align}
X^n_t=x+\int^t_0b^n_s(X^n_s)\dif s+Z_t,\quad t\in [0, 1]\label{GF2}
\end{align}
if and only if $Y^n_t=\Phi^n_t(X^n_t)$ solves the following SDE
for $t\in [0, 1]$:
\begin{align}
Y^n_t=\Phi^n_0(x)+\int^t_0a^n_s(Y^n_s)\dif s+\int^t_0\!\!\!\int_{|z|<{r_0}}
g^n_s(Y^n_{s-},z)\tilde N(\dif s,\dif z)
+\int^t_0\!\!\!\int_{|z|\geq {r_0}}g^n_s(Y^n_{s-},z) N(\dif s,\dif z),\label{Eq11}
\end{align}
where $a^n$ and $g^n$ are defined by (\ref{Lp0}) and (\ref{Lp00}).
\el

\begin{proof}
For $n\in\mN$, since $x\mapsto\Phi^n_t(x)$ and
$x\mapsto (\Phi^{n}_t)^{-1}(x)$ are smooth,
the assertion of this lemma follows from It\^o's formula.
For $n=\infty$, since we only have
$\|\nabla \Phi^\infty\|_{\infty,\gamma}<\infty$, one needs suitable mollifying technique.
This is standard and
can be found in \cite{Pr} and \cite{Zh1}. We omit the details.
\end{proof}

\bl\label{Le35}
For $n\in\mN_\infty$, let $Y^n_t(x)$ be the solution of (\ref{Eq11}) with initial value $\Phi^n_0(x)$. We have
\begin{align}
\lim_{n\to\infty}\mE\left[ \sup_{t\in[0,1]}|Y^n_t(x)-Y^\infty_t(x)|\wedge 1\right] =0.\label{Lim12}
\end{align}
Moreover, for any $p>1$, we have
\begin{align}
\sup_{n\in\mN_\infty}\sup_{x\in\mR^d}\mE\left[ \sup_{t\in[0,1]}|\nabla Y^n_t(x)|^p\right]<\infty,\label{GF1}
\end{align}
and for each $x\in\mR^d$,
\begin{align}
\lim_{n\to\infty}\mE\left[ \sup_{t\in[0,1]}|\nabla Y^n_t(x)-\nabla Y^\infty_t(x)|^p \right]=0.
\end{align}
\el
\begin{proof}
\eqref{Lim12} follows from
Lemmas \ref{L:3.2}, \ref{Le33}
and Proposition \ref{Le0} below.
In this proof, we shall drop the superscript
$``\infty''$. Notice that
\begin{align}
\nabla Y^n_t&=\nabla\Phi^n_0(x)+\int^t_0\nabla a^n_s(Y^n_s)\nabla Y^n_s\dif s+\int^t_0
\!\!\!\int_{|z|<{r_0}}\nabla_y g^n_s(Y^n_{s-},z)\nabla Y^n_{s-}\tilde N(\dif s,\dif z)\no\\
&\qquad+\int^t_0\!\!\!\int_{|z|\geq {r_0}}\nabla_yg^n_s(Y^n_{s-},z)\nabla Y^n_{s-} N(\dif s,\dif z).
\label{EA11}
\end{align}
By the Burkholder-Davis-Gundy inequality
\cite[Theorem 2.11]{Kun}
and \eqref{e:3.8}, \eqref{EK1},
we have for $p\geq 2$,
\begin{align*}
\mE\left[\sup_{s\in[0,t]}|\nabla Y^n_s|^p\right]
\preceq & |\nabla\Phi^n_0(x)|^p+\int^t_0\mE|\nabla a^n_s(Y^n_s)\nabla Y^n_s|^p\dif s
+\mE\left[\int^t_0\!\!\!\int_{|z|<{r_0}}|\nabla_yg^n_s(Y^n_{s},z)\nabla Y^n_{s}|^2 \nu(\dif z)\dif s\right]^{p/2}\\
&+\mE\left[\int^t_0\!\!\!\int_{|z|\geq {r_0}}|\nabla_yg^n_s(Y^n_{s},z)\nabla Y^n_{s}| \nu(\dif z)\dif s\right]^p
+\mE\left[ \int^t_0\!\!\!\int_{\mR^d}|\nabla_y g^n_s(Y^n_{s},z)\nabla Y^n_{s}|^p\nu(\dif z)\dif s\right] \\
\preceq & 1+\left(1+\left(\int_{|z|<{r_0}}|z|^{2\gamma}\nu(\dif z)\right)^{p/2}\right)\int^t_0\mE|\nabla Y^n_s|^p\dif s,
\end{align*}
which gives (\ref{GF1}) by Gronwall's inequality.

Next, set $U^n_t:=\nabla Y^n_t-\nabla Y_t$.
By equation \eqref{EA11}, \eqref{e:3.8}, \eqref{EK1} and \cite[Theorem 2.11]{Kun},
$$
\mE\left[\sup_{s\in[0,t]}|U^n_s|^p\right]\preceq h_n+\int^t_0\mE|U^n_s|^p\dif s,
$$
where
\begin{align*}
h_n&:= |\nabla\Phi^n_0(x)-\nabla\Phi_0(x)|^p+\int^1_0(\mE|\nabla a^n_s(Y^n_s)-\nabla a_s(Y_s)|^{2p})^{1/2}\dif s\\
&\quad +\left(\mE\left[\int^1_0\!\!\!\int_{|z|<{r_0}}|\nabla_yg^n_s(Y^n_{s},z)-\nabla_yg_s(Y_{s},z)|^2 \nu(\dif z)\dif s\right]^p\right)^{1/2}\\
&\quad+\left(\mE\left[\int^1_0\!\!\!\int_{|z\geq {r_0}}|\nabla_yg^n_s(Y^n_{s},z)-\nabla_yg_s(Y_{s},z)| \nu(\dif z)\dif s\right]^{2p}\right)^{1/2}\\
&\quad +\left(\mE\left[\int^1_0\!\!\!\int_{\mR^d}|\nabla_yg^n_s(Y^n_{s},z)-\nabla_yg_s(Y_{s},z)|^p\nu(\dif z)\dif s\right]^2\right)^{1/2}.
\end{align*}
By Gronwall's inequality,
\eqref{e:3.8}, \eqref{e:3.9}, \eqref{EK1}, \eqref{Lim4} and \eqref{Lim12},
it is easy to see that
$$
\lim_{n\to 0} \mE\left[\sup_{t\in[0,1]}|U^n_t|^p\right]\preceq \lim_{n\to 0} h_n = 0.
$$
The proof is complete.
\end{proof}

We are now in a position to give a

\begin{proof}[Proof of Theorem \ref{Main1}]
Let $a=a^\infty$ and $g=g^\infty$ be defined by (\ref{Lp0}) and (\ref{Lp00}),
respectively.  By
Lemmas \ref{L:3.2} and \ref{Le33}, we have
$$
|a_s(y)-a_s(y')|\leq C_1|y-y'|
$$
and
$$
\int_{|z|\leq r_0}|g_s(y,z)-g_s(y',z)|^2\nu(\dif z)\leq
C_2^2 |y-y'|^2\int_{|z|\leq r_0}|z|^{2\gamma}\nu(\dif z)
\stackrel{(\ref{EB1})}{\leq} C|y-y'|^2.
$$
Hence, (\ref{Eq11}) has a unique strong solution by the
classical result (cf. \cite[Theorem IV.9.1]{Ik-Wa}).
\eqref{EQ9} follows from $X_t(x)=\Phi^{-1}_t(Y_t(\Phi_0(x)))$ and \eqref{GF1}.
Moreover, let $Y_t(y)$ be the solution of SDE \eqref{Eq11} with starting point $y$.
By \eqref{e:3.10} and the choice of $r_0$ in \eqref{e:3.11},
$\{Y_t(y), t\in[0,1], y\in\mR^d\}$ defines a $C^1$-stochastic diffeomorphism flow
(cf. \cite[p.442--445]{Pr}),
so does $\{X_t(x),t\in[0,1],x\in\mR^d\}$.
Next we show that $t\mapsto \nabla X_t(x)$ is continuous. Let $X^n_t(x)$ satisfy \eqref{GF2}.
Clearly, $t\mapsto \nabla X^n_t(x)$ is continuous for each $n\in\mN$.
On the other hand, by Lemma \ref{Le35} and \eqref{Lim2}, we also have
\begin{align}
\lim_{n\to\infty}\mE\left[\sup_{t\in[0,1]} \left|\nabla X^n_t(x)-\nabla X_t(x) \right|^p\right]=0.\label{Lim13}
\end{align}
From this, we immediately obtain the desired continuity.
\end{proof}

\begin{proof}[Proof of Theorem \ref{Main2}]
First of all, we show that the right hand side of (\ref{EQ15})
is no bigger than the right hand side of (\ref{EQ16}).
By H\"older's inequality, it suffices to show that for
any $p>1$,
$$
I(t):=\mE\left[\frac{1}{S^p_t}\left|\int^t_0\nabla X_s(x)\dif
W_{S_s}\right|^p\right]\leq Ct^{- {p}/{\alpha}}.
$$
By \cite[(2.11)]{Zh3}, one has
$$
I(t)\preceq \mE\left[\frac{1}{S^p_t}\left(\int^t_0
|\nabla X_s(x)|^2\dif S_s\right)^{p/2}\right] \leq \mE\left[\frac{1}{S^{p/2}_t}\sup_{s\in[0,1]}
|\nabla X_s(x)|^p\right]
\stackrel{(\ref{EQ9})}{\preceq}\left(\mE\left[S^{-p}_t\right]\right)^{1/2}
\stackrel{(\ref{EQ18})}{\preceq} t^{-p/\alpha}.
$$
Let $b^n$ be defined as in \eqref{EQ10} and $X^n$ be the unique solution
to SDE \eqref{GF2}.
For $f\in C^1_b(\mR^d)$, by
\cite[Theorem 1.1]{Zh3} or \cite[Theorem 1.1]{Wa-Xu-Zh}, we have
$$
\nabla\mE f(X^n_t(x))=\mE\left[\frac{f(X^n_t(x))}{S_t}\int^t_0\nabla X^n_{s}(x)\dif W_{S_s}\right], \quad n\in\mN.
$$
Thus, in order to show formula (\ref{EQ15}), it suffices to show the following two relations:
\begin{equation}\label{e:new1}
\lim_{n\to\infty}\nabla\mE f(X^n_t(x))=\lim_{n\to\infty}\mE \Big[(\nabla f)(X^n_t(x))
\nabla X^n_t(x)\Big]=\mE \Big[(\nabla f)(X_t(x))\nabla X_t(x)\Big]=\nabla\mE f(X_t(x))
\end{equation}
and
\begin{align}
\lim_{n\to\infty}\mE\left[\frac{f(X^n_t(x))}{S_t}\int^t_0\nabla X^n_{s}(x)\dif W_{S_s}\right]
=\mE\left[\frac{f(X_t(x))}{S_t}\int^t_0\nabla X_{s}(x)\dif W_{S_s}\right].\label{e:new2}
\end{align}
Notice that by \eqref{Lim12} and \eqref{Lim2},
\begin{align}\label{Lim14}
\lim_{n\to\infty}\mE\left[ |X^n_t(x)-X_t(x)|\wedge 1\right]=0.
\end{align}
\eqref{e:new1} and \eqref{e:new2} follow by  \eqref{Lim13}, \eqref{Lim14} and the dominated convergence theorem.
\end{proof}

\section{Examples}

Now we give some examples for which the assumptions of Theorem \ref{Main1}
are satisfied.

\bx[\bf Subordinate Brownian motions]\label{EEU}
\rm
Let $Z_t:=W_{S_t}$, where $W$ is a $d$-dimensional  Brownian
motion with infinitesimal generator $\Delta/2$
and $S$ is a one-dimensional subordinator,
which is independent of $W_t$. Let ${\phi}(\lambda)$
be the Laplace exponent of $S$, i.e.,
$\mE \e^{-\lambda S_t}=\e^{-t{\phi}(\lambda)}$.
If for some $\alpha\in(0,2)$,
\begin{align}
{\phi}(\lambda)\geq C\lambda^{ \alpha/2},\quad \lambda\geq 1,\label{Ass}
\end{align}
then {\bf (H$^{\alpha, 1,1}_{\nu,K_0}$)} holds for some $K_0>0$.
Indeed, using the independence of $S$ and $W$, one can easily check that
for any bounded Borel function $f$ on $\mR^d$,
$$
\nabla {T}^{\nu,0}_t f(x)=\mE \left[ f(x+W_{S_t})\frac{W_{S_t}}{S_t}\right].
$$
Thus, if, for some $\beta\in (0, 1)$,
$\Lambda_x:=\sup_{y\in\mR^d}|f(x+y)-f(x)|/|y|^\beta<\infty$, then
\begin{align}
|\nabla{T}^{\nu, 0}_t f(x)|&=\left|\mE \left[\Big(f(x+W_{S_t})-f(x)\Big)\frac{W_{S_t}}{S_t}\right]\right|\no\\
&\leq \Lambda_x\mE\left[\frac{|W_{S_t}|^{1+\beta}}{S_t}\right]\leq C\Lambda_x\mE
\left[S_t^{-\frac{1-\beta}{2}}\right]\leq K_0\Lambda_x t^{(\beta-1)/\alpha},\label{KKL}
\end{align}
where the last step is due to the fact that for any
$p\in (0, 1)$,
\begin{align}
\mE S_t^{-p}&=\frac{1}{\Gamma(p)}\mE\int^\infty_0 \lambda^{p-1}
\e^{-\lambda S_t}\dif \lambda
=\frac{1}{\Gamma(p)}\int^\infty_0 \lambda^{p-1}\e^{-t{\phi}(\lambda)}
\dif \lambda\no\\
&\stackrel{(\ref{Ass})}{\leq} \frac{1}{\Gamma(p)}\left(\frac{1}{p}+\int^\infty_1
\lambda^{p-1}\e^{-C \, t\lambda^{\alpha/2}}\dif \lambda\right)\leq C
t^{-\frac{2p}{\alpha}},\ t\in(0,1].\label{EQ18}
\end{align}
The constant $C$ can be chosen to be independent of $p\in (0, 1)$ so that the constant
$K_0$ in \eqref{KKL} is independent of $\beta\in(0,1)$.
Moreover, it follows from \cite[(15)]{BGR} that
$$
\nu(\dif z)\leq  \frac{c_0 {\phi}(|z|^{-2})}{|z|^d}\dif z.
$$
Thus if there exists $\tilde\alpha\in (0, 2)$ such that
\begin{align}
{\phi}(\lambda)\leq C\lambda^{ \tilde\alpha/2} \ \
\hbox{for } \lambda\geq 1
,\label{Assub}
\end{align}
then  \eqref{EB1} is satisfied for any $\gamma \in (\tilde\alpha/2, 1]$. This implies that we need to  take
$\beta\in ( \tilde\alpha/2+1-\alpha, 1]$
in Theorem \ref{Main1}.

There are many examples of subordinate Brownian motions satisfying \eqref{Ass} and \eqref{Assub}.
One important example is the
symmetric relativistic $\alpha$-stable process
in $\mR^d$.
In this case, ${\phi}(\lambda)=(\lambda+m^{2/\alpha})^{\alpha/2}-m$ for
some $m>0$, \eqref{Ass} holds  and \eqref{Assub} is satisfied with $\tilde\alpha
=\alpha$. This implies that in this case we can take any $\beta\in (1- \alpha/2, 1]$ in Theorem \ref{Main1}.
\ex

\bx[\bf Stable-type L\'evy processes]\label{EEU1}\rm
Let ${Z}$ be a L\'evy process in $\mR^d$ whose L\'evy measure
$\nu(\dif z)= \kappa(z)\dif z$.
Assume that for some
$0<\alpha_1\leq\alpha_2<2$,
\begin{align}
c_1|z|^{-d-\alpha_1}\leq \kappa(z)\leq c_2|z|^{-d-\alpha_2}
\quad \hbox{for } |z|\leq 1.\label{Le2}
\end{align}
We call a L\'evy process satisfying the above condition of stable-type.
In this case, we can make the following decomposition for $\nu$:
$$
\nu=\nu_0+\nu_1+\nu_2
$$
with $\nu_0(\dif z):=-c_1|z|^{-d-\alpha_1} \1_{\{|z|>1\}}\dif z$ and
$$
\nu_1(\dif z):=c_1|z|^{-d-\alpha_1}\dif z,\quad
\nu_2(\dif z):=(\kappa(z)-c_1|z|^{-d-\alpha_1}) \1_{\{|z|\leq 1\}}\dif z
+ \kappa (z) \1_{\{|z|> 1\}}\dif z.
$$
By Example \ref{EEU}, ({\bf H}$^{\alpha_1, 1,1}_{\nu_1,K_0}$) holds for some $K_0>0$.
Condition \eqref{EB1} holds for any $\gamma\in ( \alpha_2/2, 1]$.
 This implies that in this case we need to take
 $\beta\in (\alpha_2/2+ 1-\alpha_1, 1]$ in Theorem \ref{Main1}.
 One particular example is the case when $\alpha_1=\alpha_2=\alpha$ and the relation in \eqref{Le2}
 is satisfied
 for all $z\in \mR^d$.
 The corresponding L\'evy process is called
 an $\alpha$-stable-like L\'evy process.
Another particular example is the case when
$\kappa (z)=0$ for $|z|>1$
 and $\alpha_1=\alpha_2=\alpha$. The corresponding
 L\'evy process is called a truncated $\alpha$-stable-like L\'evy process.
Observe that relativistic $\alpha$-stable process satisfies condition
\eqref{Le2} with $\alpha_1=\alpha_2=\alpha$.
The third particular example is the case where $\kappa (z)$ is comparable
to the L\'evy kernel of relativistic $\alpha$-stable process.
The corresponding L\'evy process can be called relativistic $\alpha$-stable-like.
 \ex

\bx[\bf Cylindrical  stable processes]\label{E:4.4}
\rm
In this example we consider a cylindrical stable process
$Z =(Z^1, \cdots, Z^k)$ in $\mR^d$, where
$Z^j$, $1\leq j\leq k$, are independent $d_j$-dimensional
rotationally symmetric
$\alpha_j$-stable processes with $\alpha_j\in (0, 2)$
and $\sum_{j=1}^kd_j=d$.
We can realize ${Z}$ as follows:
$$
Z_t=W_{S_t}:=\big(W^1_{S^1_t},\cdots,W^k_{S^k_t}\big),
$$
where $W^j$, $1\leq j\leq k$,
are independent $d_j$-dimensional standard Brownian
motions with infinitesimal generator $\Delta/2$ in $\mR^{d_j}$ and  $S^j$,
$1\leq j\leq k$,
are independent $\alpha_j/2$-stable subordinators
with $\alpha_j \in (0, 2)$ for $1\leq j\leq k$,
that are also independent of
Brownian motions $\{W^1, \dots, W^k\}$.
Define
$$
\alpha :=\min_{1\leq j\leq k}\alpha_j \quad \hbox{and} \quad \alpha_{\mathrm{max}}:=\max_{1\leq j\leq k}\alpha_j.
$$

We claim that if $\alpha \in (0, 1]$, then {\bf (H$^{\alpha, \alpha,\delta}_{\nu,K_0}$)}
holds with some $K_0>0$ and
$ \delta:=\alpha /\alpha_{\mathrm{max}}$;
and if $\alpha \in (1, 2)$, then {\bf (H$^{\alpha}_{\nu,K_0}$)} holds for some $K_0>0$.

Indeed, for $1\leq i\leq k$, let $\nabla_i=(\partial_{x_{j_i+1}}, \dots,
\partial_{x_{j_i+d_i}})$,
 where $j_i:=d_0+\cdots + d_{i-1}$ with $d_0:=0$.
As in Example \ref{EEU}, we also have the following derivative formula
for any bounded Borel function $f$ on $\mR^d$:
$$
\nabla_i {T}^{\nu, 0}_t f(x)=\mE \left[ (S^i_t)^{-1}W^i_{S^i_t}f(x+W_{S_t})\right].
$$
Suppose $\alpha\in(0,1]$. For $\beta\in[0,\alpha]$ and $x\in \mR$, if
$\Lambda_x:=\sup_{y\in\mR^d}|f(x+y)-f(x)|/|y|^\beta<\infty$,
then we have by \eqref{EQ18} that for $t\in (0, 1]$,
\begin{align*}
|\nabla_i{T}^{\nu,0}_t f(x)|
&=\left|\mE
\left[(S^i_t)^{-1}W^i_{S^i_t}\Big(f(x+W_{S_t})-f(x)\Big)\right]\right|\\
&\leq \Lambda_x\mE \left[(S^i_t)^{-1}|W^i_{S^i_t}||W_{S_t}|^\beta\right]\\
&\leq \Lambda_x\left(\mE\left[(S^i_t)^{-1}|W^i_{S^i_t}|^{1+\beta}\right]
+\mE\left [(S^i_t)^{-1}|W^i_{S^i_t}|\right]\sum_{j\not= i}\mE\Big[ |W^j_{S^j_t}|^\beta \Big]\right)\\
&\leq C\Lambda_x\left(
t^{(\beta-1)/\alpha_i}+t^{-1/\alpha_i} \sum_{j\not= i}t^{\beta /\alpha_j}
\, \mE\Big[ |W^j_{S^j_1}|^\beta \Big]\right) \\
&\leq K_0\Lambda_x\left(
t^{{(\beta-1)}/{\alpha}}
+t^{\beta/\alpha_{\mathrm{max}}-1/\alpha }\right)
\qquad \hbox{(since $\beta<\alpha$)}\\
&\leq K_0\Lambda_x
t^{{(\alpha\beta/\alpha_{\mathrm{max}}-1)}/{\alpha}}
=K_0\Lambda_x t^{(\delta\beta-1)/\alpha};
\end{align*}
that is, {\bf (H$^{\alpha, \alpha,\delta}_{\nu,K_0}$)} holds.
If $\alpha\in(1,2)$, then we have by \eqref{EQ18} that for $t\in (0, 1]$,
\begin{align*}
 |\nabla_i{T}^{\nu, 0}_t f(x)|
\leq\|f\|_\infty\mE
\left[  |W^i_{S^i_t}|/ (S^i_t) \right]
\preceq \|f\|_\infty\mE \left[(S^i_t)^{-1/2}\right]
\preceq \|f\|_\infty t^{-1/\alpha_i}\leq K_0\|f\|_\infty t^{-1/\alpha}.
\end{align*}
Thus in this case, {\bf (H$^{\alpha}_{\nu,K_0}$)} holds. The claim is now verified.

It is not difficult to see by using the property of the rotationally symmetric $\alpha_j$-stable process $W^j_{S_j}$ that the parameter $\alpha$ in
the now verified property {\bf (H$^{\alpha, \alpha,\delta}_{\nu,K_0}$)} and {\bf (H$^{\alpha}_{\nu,K_0}$)} is best possible. For example, it can be shown that
when $\alpha \in (1, 2)$, property {\bf (H$^{\alpha^*}_{\nu,K_0}$)}
fails for any $\alpha^*>\alpha$.

Note that \eqref{EB1} holds for any $\gamma >\alpha_{\mathrm{max}}/2$.
For Theorem \ref{Main1} to be valid,
the following constraint  needs to be satisfied:
 $$
1 \geq \beta > \alpha_{\mathrm{max}}/2 + \alpha_{\mathrm{max}} (1-\alpha)/\alpha
\ \mbox{ if } \alpha\leq 1,
 \quad \mbox{and} \quad
\alpha_{\mathrm{max}}<2\alpha\ \mbox{ if } \alpha>1 .
 $$
Clearly, when $\alpha>1$, the condition $\alpha_{\mathrm{max}}<2 \alpha$
is automatically satisfied.
Consequently, in this case for Theorem \ref{Main1} to be applicable,
we need $\alpha_i$'s to satisfy
\begin{equation}\label{e:1.14}
\hbox{either} \quad \alpha>1 \quad \hbox{or} \quad \alpha \in (0, 1]
\ \hbox{ and } \ \alpha_{\mathrm{max}}<2\alpha^2/(2-\alpha),
 \end{equation}
and take
$$
\beta \in (\beta_0, 1]\
\hbox{ with }\beta_0:=\alpha_{\mathrm{max}}/2+(\alpha_{\mathrm{max}}/\alpha
\1_{\{\alpha\leq 1\}} + \1_{\{\alpha>1\}})(1-\alpha).
$$

Condition \eqref{e:1.14} implies that $\alpha >2/3$.
An open question is   whether
constraint \eqref{e:1.14} can be dropped.
It boils down to the question whether
{\bf (H$^{\alpha, 1,1}_{\nu,K_0}$)} holds for any cylindrical stable process.

This example can be extended in two directions. First,
as in Example \ref{EEU}, we can consider more general subordinators
$\{S^1, \dots, S^k\}$.
Second,
as in Example \ref{EEU1}, we can consider more general L\'evy process,
whose L\'evy measure is bounded by the L\'evy measure of the cylindrical
$\alpha$-stable process $W_{S}$ (or, more generally, the cylindrical
subordinate Brownian motion) from below.
\ex

\bigskip

\noindent{\bf Proof of Corollary \ref{C:1.3}.}
It follows from Examples \ref{EEU}, \ref{EEU1} and \ref{E:4.4}. \qed

\section{Appendix}

In this appendix, we prove the continuous dependence of solutions to SDEs with jumps
with respect to the initial values and coefficients.
\bp\label{Le0}
Fix $r>0$. Let $a^n, g^n, n\in\mN_\infty$ be two families of
uniformly Lipschitz continuous functions in the sense that
for some $C>0$, and all $n\in\mN_\infty$ and $t\in[0,1], x,y,z\in\mR^d$,
\begin{align}
|a^n_t(x)-a^n_t(y)|\leq C|x-y|,\ \ |g^n_t(x,z)-g^n_t(y,z)|\leq C|x-y| h(z),\label{EQ12}
\end{align}
where $\int_{|z|\leq r}|h(z)|^2\nu(\dif z)<\infty$.
Suppose that for each $t\in[0,1]$ and $x,z\in\mR^d$,
\begin{align}
\lim_{n\to\infty}a^n_t(x)=a^\infty_t(x), \quad
\lim_{n\to\infty}g^n_t(x,z)=g^\infty_t(x,z)\label{EQ13}
\end{align}
and
\begin{align}
\sup_{n\in\mN_\infty}\sup_{t\in[0,1]}\sup_{x\in\mR^d}\left(\frac{|a^n_t(x)|}{1+|x|}+\sup_{0<|z|\leq r}\frac{|g^n_t(x,z)|}{|z|}\right)<\infty.\label{EQ133}
\end{align}
For $n\in\mN_\infty$, let $Y^n_t$ be the solution to the following SDE
$$
Y^n_t=\xi_n+\int^t_0a^n_s(Y^n_s)\dif s+\int^t_0\!\!\!
\int_{|z|\leq r}
g^n_s(Y^n_{s-},z)\tilde N(\dif s,\dif z)
+\int^t_0\!\!\!
\int_{|z|> r}
g^n_s(Y^n_{s-},z) N(\dif s,\dif z).
$$
Assume that $\xi_n$ converges to $\xi_\infty$ in probability as $n\to\infty$. Then we have
\begin{align}
\lim_{n\to\infty}\mE\left(\sup_{t\in[0,1]}|Y^n_t-Y^\infty_t|\wedge 1\right)=0,\label{EQ14}
\end{align}
which implies that $Y^n_t$ converges to $Y^\infty_t$ in probability.
\ep

We begin with the following lemma.
\bl\label{Le52}
There is a nonnegative smooth function $f$ on $\mR^d$ with the following properties:
\begin{align}\label{GF4}
f(x)=|x|^2 \ \hbox{ if }  |x|\leq 1,\quad  f(x)=2 \ \hbox{ if } |x|\geq 2,
 \quad \hbox{and} \quad |\nabla f|+|\nabla^2 f|
\leq C_1 \1_{\{|x|\leq 2\}},
\end{align}
for some constant $C_1>0$, and that for any constant $C_2>0$,
there exists a constant $C_3>0$ such that
for all  $\delta>0$, $r\in[0,1]$ and
$|y|\leq C_2((|x|+\delta)\wedge 1)$,
\begin{align}
|y||\nabla f(x+ry)|\leq C_3(f(x)+\delta), \quad |y|^2|\nabla^2f(x+ry)|\leq C_3(f(x)+\delta^2).\label{GF3}
\end{align}
\el
\begin{proof}
Let $\phi$ be an increasing smooth function on $(0,\infty)$ with $\phi(r)=r$ for $r\leq 1$
and $\phi(r)=2$ for $r\geq 4$. Let $f(x):=\phi(|x|^2)$. It is easy to check that
$f$ has the desired properties.
\end{proof}

We also need the following key lemma.
\bl\label{Le53}
Let $\tau_1, \tau_2$ be two stopping times with $0\leq\tau_1\leq \tau_2\leq 1$.
In the setup of  Proposition \ref{Le0},
let $Y^n$ solve the following SDE on $[\tau_1,\tau_2]$:
$$
Y^n_t=Y^n_{\tau_1}+\int^t_{\tau_1}a^n_s(Y^n_s)\dif s+\int^t_{\tau_1}\!\int_{|z|\leq r}
g^n_s(Y^n_{s-},z)\tilde N(\dif s,\dif z).
$$
Assume that $Y^n_{\tau_1}$ converges to $Y^\infty_{\tau_1}$ in probability, then we have
$$
\lim_{n\to\infty}\mE\left[ \sup_{t\in[\tau_1,\tau_2]}|Y^n_t-Y^\infty_t|\wedge 1\right] =0.
$$
\el

\begin{proof}
In this proof we will drop the superscript ``$\infty$'' and write
$$
U^n_s:=Y^n_s-Y_s,\quad A^n_s:=a^n_s(Y^n_s)-a_s(Y_s), \quad
\Gamma^n_s(z):=g^n_s(Y^n_{s-},z)-g_s(Y_{s-},z).
$$
Let $f$ be as in Lemma \ref{Le52}. By It\^o's formula, we have
\begin{align*}
f(U^n_t)&=f(U^n_{\tau_1})+\int^t_{\tau_1}\<A^n_s,\nabla f(U^n_s)\>\dif s
+\int^t_{\tau_1}\!\int_{|z|\leq r}[f(U^n_{s-}+\Gamma^n_s(z))-f(U^n_{s-})]\tilde N(\dif s,\dif z)\\
&\quad+\int^t_{\tau_1}\!\int_{|z|\leq r}[f(U^n_{s-}+\Gamma^n_s(z))-f(U^n_{s-})-\Gamma^n_s(z)\cdot\nabla f(U^n_{s-})]
\nu(\dif z)\dif s.
\end{align*}
For $R>0$, define a stopping time
$$
\tau_R:=\inf\Big\{t\geq\tau_1: |Y_s|>R\Big\}\wedge\tau_2.
$$
For any $T\in[0,1]$,
by the Burkholder-Davis-Gundy inequality \cite[Theorem 2.11]{Kun}, we have
\begin{eqnarray*}
&&\mE\left(\sup_{t\in[\tau_1,T\wedge\tau_R]}|f(U^n_t)|^2\right)\\
&&\preceq
\mE|f(U^n_{\tau_1})|^2+\mE\int^{T\wedge\tau_R}_{\tau_1}|\<A^n_s,\nabla f(U^n_s)\>|^2\dif s\\
&& \quad +\mE\int^{T\wedge\tau_R}_{\tau_1}\!\!\!\int_{|z|\leq r}|f(U^n_{s-}+\Gamma^n_s(z))-f(U^n_{s-})|^2\nu(\dif z)\dif s\\
&&\quad +\mE\int^{T\wedge\tau_R}_{\tau_1}\!\left|\int_{|z|\leq r}[f(U^n_{s-}+\Gamma^n_s(z))-f(U^n_{s-})-
\Gamma^n_s(z)\cdot\nabla f(U^n_{s-})]
\nu(\dif z)\right|^2\dif s\\
&&=: \mE|f(U^n_{\tau_1})|^2+I_1^n+I^n_2+I^n_3.
\end{eqnarray*}
For $I^n_1$,
by \eqref{EQ12} and  \eqref{GF4},
we have
$$
I^n_1\preceq \mE\int^{T\wedge\tau_R}_{\tau_1}|a^n_s(Y_s)-a_s(Y_s)|^2\dif s+\mE
\int^{T\wedge\tau_R}_{\tau_1}|f(U^n_s)|^2\dif s.
$$
For $I^n_2$ and $I^n_3$,
noticing that  by \eqref{EQ12} and \eqref{EQ133},
$$
|\Gamma^n_s(z)|\leq C((|U^n_{s-}|+|g^n_s(Y_{s-},z)-g_s(Y_{s-},z)|)\wedge 1),
\quad |z|\leq r,
$$
by \eqref{GF3},
we have
\begin{align*}
I^n_2&\preceq \mE\int^{T\wedge\tau_R}_{\tau_1}\!\!\!\int_{|z|\leq r}|\Gamma^n_s(z)|^2\left(\int^1_0|\nabla f(U^n_{s-}+r\Gamma^n_s(z))|^2\dif r\right)\nu(\dif z)\dif s\\
&\preceq\mE\int^{T\wedge\tau_R}_{\tau_1}\!\!\!\int_{|z|\leq r}|g^n_s(Y_s,z)-g_s(Y_s,z)|^2\nu(\dif z)\dif s
+\mE\int^{T\wedge\tau_R}_{\tau_1}|f(U^n_s)|^2\dif s,
\end{align*}
and
\begin{align*}
I^n_3&\preceq \mE\int^{T\wedge\tau_R}_{\tau_1}\left(\int_{|z|\leq r}|\Gamma^n_s(z)|^2
\left(\int^1_0\!\!\!\int^1_0|\nabla^2 f(U^n_{s-}+rr'\Gamma^n_s(z))|
\dif r\dif r'\right)\nu(\dif z)\right)^2\dif s\\
&\preceq\mE\int^{T\wedge\tau_R}_{\tau_1}\left(\int_{|z|\leq r}|g^n_s(Y_s,z)-g_s(Y_s,z)|^2
\nu(\dif z)\right)^2\dif s
+\mE\int^{T\wedge\tau_R}_{\tau_1}|f(U^n_s)|^2\dif s.
\end{align*}
Combining the above calculations, we obtain
$$
\mE\left[ \sup_{t\in[\tau_1,T\wedge\tau_R]}|f(U^n_t)|^2\right] \preceq
h_n+\mE\int^{T\wedge\tau_R}_{\tau_1}|f(U^n_s)|^2\dif s,
$$
where
\begin{align*}
h_n&:=
\mE \left[ |f(U^n_{\tau_1})|^2 \right]
+\mE\int^{T\wedge\tau_R}_{\tau_1}|a^n_s(Y_s)-a_s(Y_s)|^2\dif s \\
&\quad +\mE\int^{T\wedge\tau_R}_{\tau_1}\!\!\!\int_{|z|\leq r}|g^n_s(Y_s,z)-g_s(Y_s,z)|^2\nu(\dif z)\dif s\\
&\quad +\mE\int^{T\wedge\tau_R}_{\tau_1}\left(\int_{|z|\leq r}|g^n_s(Y_s,z)-g_s(Y_s,z)|^2\nu(\dif z)\right)^2\dif s.
\end{align*}
By Gronwall's inequality, \eqref{EQ13}, \eqref{EQ133}  and the dominated convergence theorem, we have for each $R>0$,
$$
\lim_{n\to 0} \mE\left[ \sup_{t\in[\tau_1,\tau_R]}|f(U^n_t)|^2\right] \leq
\lim_{n\to 0} Ch_n =0.
$$
In particular,
$$
\lim_{n\to\infty}\mE\left[ \sup_{t\in[\tau_1,\tau_R]}|Y^n_t-Y^\infty_t|^4\wedge 1\right] =0,
$$
which together with $\lim_{R\to\infty} \mP(\tau_R<\tau_2)=0$ gives the desired limit.
\end{proof}

Now we give
\begin{proof}[Proof of Proposition \ref{Le0}]
Let $\tau_1:=0$ and for $m\in\mN$, define recursively
$$
\tau_{m+1}:=\inf\{t>\tau_m: |Z_s-Z_{s-}| > r\}\wedge 1.
$$
Since $Z$ only has finite many jumps greater than $r$ before time $1$, we have
$\lim_{m\to\infty}\tau_m=1$. Clearly, for $t\in(\tau_{m},\tau_{m+1}]$, $Y^n_t$ satisfies
$$
Y^n_t=Y^n_{\tau_m}+\int^t_{\tau_{m}}a^n_s(Y^n_s)\dif s+\int^t_{\tau_{m}}\!
\int_{|z|\leq r}
g^n_s(Y^n_{s-},z)\tilde N(\dif s,\dif z),
$$
where $Y^n_{\tau_m}:=Y^n_{\tau_{m}-}+g^n_{\tau_m}(Y^n_{\tau_{m}-}, Z_{\tau_{m}}-Z_{\tau_m-})$.
Since $\xi_n\to \xi_\infty$ in probability as $n\to\infty$, by  Lemma \ref{Le53}
and  induction, we have for each $m\in\mN$,
$$
\lim_{n\to\infty}\mE\left[ \sup_{t\in[\tau_m,\tau_{m+1}]}
|Y^n_t-Y^\infty_t|\wedge 1\right]=0,
$$
which gives the desired limit by
noticing that for any $m_0\in\mN$,
$$
\mE\left[ \sup_{t\in[0,1]}|Y^n_t-Y^\infty_t|\wedge 1\right]
\leq\sum_{m=1}^{m_0}\mE\left[ \sup_{t\in[\tau_m,\tau_{m+1}]}|Y^n_t-Y^\infty_t|\wedge 1\right] +
\mP(\tau_{m_0+1}<1),
$$
and $\lim_{m_0\to\infty} \mP(\tau_{m_0+1}<1)=0$.
\end{proof}

\bigskip

{\bf Acknowledgement.} We thank Zenghu Li, Feng-Yu Wang and Jian Wang for helpful discussions.

\bigskip


\begin{thebibliography}{999}

\bibitem{Ap}Applebaum, D.:
{\it L\'evy Processes and Stochastic Calculus}.
Cambridge Studies in Advance Mathematics 93, Cambridge University Press, 2004.

\bibitem{BGR} Bogdan, K., Grzywny, T. and Ryznar, M.:
Density and tails of unimodal convolution semigroiups.
{\it J. Funct. Anal.} {\bf 266} (2014), 3543--3571.

\bibitem{Ch-Zh1}Chen, Z.-Q. and Zhang, X.:
Heat kernels and analyticity of
non-symmetric jump diffusion semigroups. Preprint.  arXiv:1306.5015v2 [math.AP]

\bibitem{Fe-Fl}Fedrizzi, E. and Flandoli F.:
Pathwise uniqueness and
continuous dependence for SDEs with nonregular drift.
Arxiv:1004.3485v1.

\bibitem{Fl-Gu-Pr}Flandoli, F., Gubinelli, M. and Priola, E.:
Well-posedness
of the transport equation by stochastic perturbation. {\it Inven. Math. \bf 180} (2010), 1-53.

\bibitem{Ik-Wa}Ikeda, N., Watanabe, S.:
{\it Stochastic Differential Equations and Diffusion Processes},
2nd ed.. North-Holland/Kodanska,
Amsterdam/Tokyo, 1989.

\bibitem{Kr}Krylov, N.V.:
{\it Lectures on Elliptic and Parabolic Equations in H\"older Spaces}.
Graduate Studies in Mathematics, Vol. 12, AMS, 1996.


\bibitem{Kr-Ro}Krylov, N.V. and R\"ockner, M.:
Strong solutions of
stochastic equations with singular time dependent drift.
{\it Probab. Theory Relat. Fields \bf 131} (2005), 154-196.

\bibitem{Kun} Kunita, H.:
Stochastic differntial equations based on L\'evy processes
and stochastic flows of diffeomorphisms. In \emph{Real
and Stochastic Analysis}, 305--373. Birkh\"auser, Boston, 2004.

\bibitem{Pr}Priola, E.:
Pathwise uniqueness for singular SDEs driven by
stable processes. {\it Osaka Journal of Mathematics}, {\bf 49} (2012), 421-447.

\bibitem{Pr1}Priola, E.:
Stochastic flow for SDEs with jumps and irregular
drift term. arXiv:1405.2575v1.

\bibitem{Pro}Protter, P.:
{\it Stochastic Integration and Differential
Equations}. 2nd ed., Springer-Verlag, Berlin, 2004.

\bibitem{Si2}Silvestre, L.:
On the differentiability of the solution
to an equation with drift and fractional diffusion.
{\it Indiana Univ. Math. J. \bf 61} (2012), 557-584.

\bibitem{St}Stein, E.M.:
{\it Singular Integrals and
Differentiability Properties of Functions}. Princeton, N.J.,
Princeton University Press,  1970.

\bibitem{Ta-Ts-Wa}Tanaka, H., Tsuchiya, M. and Watanabe, S.:
Perturbation of drift-type for L\'evy processes.
{\it J. Math. Kyoto Univ. \bf 14} (1974), 73-92.


\bibitem{Ve}Veretennikov, A. Ju.:
On the strong solutions of
stochastic differential equations. {\it Theory Probab. Appl. \bf
24} (1979), 354-366.

\bibitem{Wa-Xu-Zh}Wang, F., Xu, L. and Zhang, X.: Gradient
estimates for SDEs driven by multiplicative L\'evy noise.
arXiv:1301.4528.

\bibitem{Zh0}Zhang, X.: Stochastic homeomorphism flows of
SDEs with singular drifts and Sobolev diffusion coefficients.
{\it Electron. J. Probab. \bf 16} (2011), 1096-1116.

\bibitem{Zh2}Zhang, X.: Stochastic functional differential
equations driven by
L\'evy processes and quasi-linear partial integro-differential
equations. {\it Ann. Appl. Aprob. \bf 22} (2012), 2505-2538.

\bibitem{Zh1}Zhang, X.: Stochastic differential equations with
Sobolev drifts and driven by $\alpha$-stable processes.
{\it Ann. Inst. H. Poincare Probab. Statist. \bf 49} (2013), 1057-1079.

\bibitem{Zh3}Zhang, X.: Derivative formula and gradient estimate for SDEs
driven by $\alpha$-stable processes. {\it Stoch. Proc. Appl.}
{\bf 123} (2013), 1213-1228.

\bibitem{Zv}Zvonkin, A.K.: A transformation of the phase
space of a diffusion process
that removes the drift.  {\it Mat. Sbornik},
{\bf 93 (135)} (1974), 129-149.

\end{thebibliography}
\end{document}